\documentclass[10pt,oneside]{article}
\usepackage{latexsym,amssymb,amsmath}
\usepackage{graphicx}
\usepackage{color}
\usepackage[emtex,matrix,arrow]{xy}

\newtheorem{thm}{Theorem}
\newtheorem{cor}[thm]{Corollary}
\newtheorem{ex}[thm]{Example}

\newtheorem{lem}[thm]{Lemma}
\newtheorem{prop}[thm]{Proposition}
\newtheorem{rem}[thm]{Remark}
\numberwithin{equation}{section}

\newcommand{\kb}{\mathbb K}

\newcommand{\espai}{\vspace{0.3 truecm}}

\newcommand{\To}{\rightarrow}

\newcommand{\Aa}{\mathcal{A}}

\newcommand{\Gg}{\mathcal{G}}

\newcommand{\Cc}{\mathcal{C}}

\newcommand{\Ee}{\mathcal{E}}

\newcommand{\F}{\mathbb{F}}
\newcommand{\Vv}{\mathcal{V}}

\newcommand{\G}{\mathbb{G}}

\newcommand{\Cgg}{\mathbf{C}}

\newcommand{\ev}{\mathcal{V}ect}

\newcommand{\qed}{\hspace*{\fill}$\Box$  \ifmmode \else
    \par\addvspace\topsep\fi}

\newenvironment {proof}{\par\addvspace\topsep\noindent{\it Proof.}
   \ignorespaces }{\qed}

\begin{document}

\title{The 2-group of symmetries of a split chain complex} 
\author{Josep Elgueta \thanks{This work was partially supported by the Generalitat de Catalunya (Project: 2009 SGR 1284) and the Ministerio de Educaci\'on y Ciencia of Spain (Project: MTM2009-14163-
C02-02).} \\ Dept. Matem\`atica Aplicada II \\ Universitat Polit\`ecnica de Catalunya}
\date{}

\maketitle

\abstract{
We explicitly compute the 2-group of self-equivalences and (homotopy classes of) chain homotopies between them for any {\it split} chain complex $A_{\bullet}$ in an arbitrary $\kb$-linear abelian category ($\kb$ any commutative ring with unit). In particular, it is shown that it is a {\it split} 2-group whose equivalence class depends only on the homology of $A_{\bullet}$, and that it is equivalent to the trivial 2-group when $A_\bullet$ is a split exact sequence. This provides a description of the {\it general linear 2-group} of a Baez and Crans 2-vector space over an arbitrary field $\mathbb{F}$ and of its generalization to chain complexes of vector spaces of arbitrary length.}

\section{Introduction}

In the last years, several clues have appeared suggesting that the basic set-up of representation theory should be widened. One way of doing this consists of representing groups not as symmetries of objects in a category, but as symmetries of objects in a {\it 2-category}.

Roughly, a 2-category is the same thing as a category, except that we also have morphisms between the morphisms (or 2-morphisms), and two associated ways of composing them, called {\it vertical} and {\it horizontal} compositions, depending on the ``dimension'' of the common cell (see Figure~1).

\begin{figure}[htbp!]
\centering
\input{composicions.pstex_t}
\caption{}
\label{figura_llei_intercanvi}
\end{figure}

\noindent
Thus objects $X$ in a 2-category $\mathbf{C}$ not only have symmetries, but also symmetries between the symmetries, so that we have a {\it groupoid} of symmetries for them. Moreover, this groupoid comes equipped with a natural {\it monoidal structure}, i.e. an associative and unital (up to isomorphism) product, given by the composition of morphisms and the horizontal composition of 2-morphisms. Hence we have a sort of categorified group. It will be called the {\it 2-group of symmetries} of $X$. More generally, by a {\it 2-group} it is meant any groupoid $\Gg$ equipped with a product functor $\Gg\times\Gg\To\Gg$ which is associative and unital up to given (coherent) isomorphisms and such that each object has a (possibly weak) inverse. When $\Gg$ is discrete (only identity morphisms) this reduces to the usual notion of group.

The paradigmatic example of 2-category $\mathbf{C}$ is the 2-category $\mathbf{Cat}$ with all (small) categories as objects, functors between them as morphisms, and natural transformations between these as 2-morphisms. Then the 2-group of symmetries of a given category $\Cc$ is the groupoid having as objects all self-equivalences of $\Cc$ and as morphisms all natural isomorphisms between these. The product functor is given by composition of self-equivalences and the usual horizontal composition of natural isommorphisms, and the unit is the identity of $\Cc$. If $\Cc$ is the discrete category defined by some set $S$, we recover the group of automorphisms of $S$.

If we take this idea seriously, the first problem we should face is deciding on what 2-category we wish to represent groups and more generally, 2-groups. This naturally leads to the search of an analog to the category of vector spaces in this new setting. Let us call it the 2-category of {\it 2-vector spaces} (over a given field $\mathbb{F}$).

There are various more or less natural canditates to the notion of 2-vector space. Actually, this is typical of the categorification of any given algebraic structure. The most popular one is that introduced by Kapranov and Voevodsky \cite{KV94}. These authors get to the notion of 2-vector space by categorifying the usual notion of vector space over a field $\mathbb{F}$. They take as analog of $\mathbb{F}$ the (semiring) category $\ev_{\mathbb{F}}$ of finite dimensional vector spaces over $\mathbb{F}$ and consider symmetric monoidal categories (categorical analogs of the abelian groups) on which $\ev_\mathbb{F}$ acts. These are called $\ev_\mathbb{F}$-{\it module categories}. In fact, they restrict to the ``free'' such objects, i.e. those of the form $\ev_\mathbb{F}^n$ for some $n\geq 0$. The representation theory of (2-)groups in these 2-vector spaces is studied in \cite{jE4}, \cite{jE6}. 

Kapranov and Voevodsky 2-vector spaces can be generalized in various ways (see for instance \cite{jE5}). There is, however, another sensible way of defining 2-vector space. It is based on the fact that categories are nothing but simplicial sets of a particular kind. Thus, as Grothendieck first pointed out, any category is completely described (up to {\it isomorphism}) by its {\it nerve}. This is the simplicial set having as $n$-simplices all paths of morphisms of length $n$, and as face and degeneracy maps those given by composition and insertion of identity morphisms, respectively. In fact, the simplicial sets arising as nerves of a category admit a very neat characterization as those satisfying the so called {\it Segal condition} or {\it nerve condition}~\footnote{For any $k\geq 2$, there is a canonical map between the set of paths of 1-cells of length $k$ in the simplicial set and its set of $k$-cells (see for instance \cite{LP08}). These are the so called {\it Segal maps}. The nerve condition is that all these Segal maps must be bijections.}. Anyway, what is relevant for our purposes is that we can go from sets to categories by just going from sets to (some kind of) {\it simplicial} sets. This suggests that a sensible notion of 2-vector space should be given by (some kind of) simplicial vector spaces. Now, according to the Dold-Kan correspondence, simplicial objects in any abelian category $\Aa$ are equivalent to positive chain complexes of objects in $\Aa$ (see \cite{cW94}). Therefore we are led to think of the chain complexes of vector spaces (or at least, some special type of them) as a suitable notion of 2-vector space, and to take as 2-category on which to represent 2-groups a suitably defined 2-category $\mathbf{Ch}(\ev_\mathbb{F})$ with the chain complexes of vector spaces as objects.

In fact, in \cite{BC03} Baez and Crans defined a 2-vector space as a category in the category of vector spaces~\footnote{This is sometimes called {\it internal categorification} and it provides a general method of categorifying certain objects. The idea is that a structure defined on the category of sets can be categorified by taking the same structure in the category of categories. However, this only gives {\it strict} categorical versions.} and proved that this amounts to a chain complex of length 1. As argued before, however, it looks reasonable to consider as 2-vector spaces chain complexes of vector spaces of arbitrary length.

In this paper we adopt this point of view, and we undertake the task of computing explicitly the corresponding {\it simplicial general linear 2-groups}. This should be viewed as a preliminary step toward a {\it simplicial} representation theory for 2-groups. Indeed, such a representation will be a morphism of 2-groups into some of these simplicial general linear 2-groups.

Apart from the {\it discrete} and {\it one-object} 2-groups, defined below and which just amount to groups and abelian groups respectively, the next simplest kind of 2-groups are the so called {\it split}. These are the 2-groups whose underlying monoidal groupoid has a skeleton which is a {\it strict} monoidal subgroupoid (i.e. the associativity and unit constraints for the restricted product functor are identities). As a matter of fact, the representation theory of 2-groups becomes much easier when the 2-category on which we represent them is such that the 2-groups of symmetries of its objects are split. The 2-category of Kapranov and Voevodsky 2-vector spaces and its generalization introduced in \cite{jE5} is of this kind. In this work we show that this is also the case for the 2-category $\mathbf{Ch}(\ev_\mathbb{F})$ and more generally, for the 2-category $\mathbf{Ch}(\Aa)$ of chain complexes in any semisimple abelian category $\Aa$.

\espai
The outline of the paper is as follows. Section 2 is a quick review on 2-groups. Apart from the basic definitions, Sinh's classification theorem is reviewed in some detail. In Section 3 we recall the notion of (split) chain complex of objects in an arbitrary abelian category $\Aa$ and describe the corresponding 2-category $\mathbf{Ch}(\Aa)$ in elementary terms. Finally, in Section 4 the 2-group of symmetries of an arbitrary split object $A_{\bullet}$ of $\mathbf{Ch}(\Aa)$ is computed (cf. Theorem~\ref{teorema_principal}), and an explicit equivalence between this 2-group and its equivalent split version is described.

\espai
We assume the reader is familiar with the notion of (weak) monoidal category; see for instance \cite{sM98}. The unit object, associator and left and right unit isomorphisms are denoted by $I$, $a$, $l$ and $r$, respectively. The corresponding many objects version, i.e. the notion of (weak) 2-category is recalled in the Appendix.

All over the paper, abelian categories are assumed to be over an arbitrary commutative unital ring $\kb$. Thus all hom-sets are $\kb$-modules and all composition maps are $\kb$-bilinear. The reader may think of the case $\kb=\mathbb{Z}$.

\section{Review on 2-groups}
\label{2grups}

\subsection{Basics on 2-groups}

Roughly speaking, a {\it 2-group} (also called {\it categorical group} or {\it gr-category}) is a category equipped with a (suitably weakened) group structure. More precisely, it is a monoidal category
$$\mathbb{G}=(\Gg,\otimes,I,a,l,r)
$$
whose underlying category $\Gg$ is a groupoid and such that for each object $x$ there exists a weak inverse, i.e. an object $x^*$ such that
$$
x\otimes x^*\cong I\cong x^*\otimes x.
$$
When the monoidal category is strict (i.e. the associator $a$ and the left and right unit constraints $l$ and $r$ are identities) and all inverses $x^*$ are strict (i.e. $x\otimes x^*=I=x^*\otimes x$), $\mathbb{G}$ is called a {\it strict} 2-group. For more details, see \cite{BL03}.

The simplest examples of 2-groups are the groups themselves thought of as discrete categories. For any group $G$ we shall denote by $G[0]$ the corresponding {\it discrete 2-group}. If the group is abelian, it can also be viewed as a {\it one-object 2-group}. In this case, the group plays the role of the group of automorphisms of the unique object. For any abelian group $A$, we shall denote by $A[1]$ the corresponding one-object 2-group. In both cases, the tensor product is given by the group law, and both are examples of strict 2-groups.

More generally, for any group $G$ and any $G$-module $A$ we have a strict 2-group defined as follows (the previous cases correspond to $A$ or $G$ trivial). Its set of objects is $G$, its set of morphisms is $A\times G$, with $(a,g)$ an automorphism of $g$, and composition and tensor product are given by
\begin{align*}
(a',g)\circ(a,g)&=(a'+a,g),
\\ g\otimes g'&=gg', \\ (a,g)\otimes(a',g')&=(a+g\lhd a',gg'),
\end{align*}
where $\lhd:G\times A\To A$ denotes the action of $G$ on $A$. This is just a special case of the general notion of semidirect product for 2-groups, in this case between $G[0]$ and $A[1]$ (see \cite{GI01}). The 2-group so defined will be denoted by $A[1]\rtimes G[0]$, and its underlying groupoid by $\Gg_{A,G}$. Notice that it is a strict 2-group and that $\Gg_{A,G}$ is skeletal (isomorphic objects are equal). Any 2-group of this kind or equivalent (in a sense made precise below) to one of this kind will be called a {\it split 2-group}.

2-groups arise naturally as symmetries of objects in a 2-category (for the precise definition of a 2-category, cf. Appendix). Thus for any 2-category $\mathbf{C}$ and any object $X$ of $\mathbf{C}$ the groupoid $\mathcal{E}quiv(X)$ of self-equivalences of $X$ and 2-isomorphisms between these has a canonical structure of a 2-group. The tensor product is given by the composition of self-equivalences and the horizontal composition of 2-morphisms. We shall denote the 2-group so defined by $\mathbb{E}quiv(X)$. The unit object is $id_X$. Notice that the underlying monoidal groupoid is strict when $\mathbf{C}$ is strict. However, even in this case $\mathbb{E}quiv(X)$ may be non-strict because there may exist objects with no strict inverse, i.e. self-equivalences of $X$ which are not automorphisms.

2-groups are the objects of a 2-category $\mathbf{2Grp}$ whose 1-morphisms are the monoidal functors between the underlying monoidal groupoids. Thus, for any 2-groups $\mathbb{G}$ and $\mathbb{G}'$, a 1-morphism from $\G$ to $\G'$ is given by a pair $\F=(F,\mu)$ with $F:\Gg\To\Gg'$ a functor and $\mu$ a collection of natural isomorphisms (the {\it monoidal structure})
$$
\mu_{x,y}:F(x\otimes y)\stackrel{\cong}{\To}F(x)\otimes' F(y)
$$
labelled by pairs of objects of $\Gg$ and such that the diagram
\begin{equation} \label{axioma_coherencia}
\xymatrix{
F((x\otimes y)\otimes z)\ar[rr]^{F(a_{x,y,z})}\ar[d]_{\mu_{x\otimes y,z}} && F(x\otimes(y\otimes z))\ar[d]^{\mu_{x,y\otimes z}} \\ F(x\otimes y)\otimes' F(z)\ar[d]_{\mu_{x,y}\otimes' id_{F(z)}} && F(x)\otimes' F(y\otimes z)\ar[d]^{id_{F(x)}\otimes'\mu_{y,z}} \\ (F(x)\otimes' F(y))\otimes' F(z)\ar[rr]_{a'_{F(x),F(y),F(z)}} && F(x)\otimes'(F(y)\otimes' F(z)) 
}
\end{equation}
commutes for all objects $x,y,z\in\G$. As in the case of groups, the unit object as well as the inverses are automatically preserved, at least up to isomorphism. Notice that we do not mention explicitly the unit isomorphism $F(I)\cong I'$ usually required in the definition of a monoidal functor. This is because it is uniquely determined by the $\mu$'s when the monoidal functor is between 2-groups instead of arbitrary monoidal categories. As it concerns 2-morphisms, they are given by the so called {\it monoidal natural transformations} between these monoidal functors (see \cite{sM98} or \cite{BL03} for the precise definition).

Since 2-groups are the objects of a 2-category it makes sense to speak of {\it equivalent} 2-groups, i.e. 2-groups $\mathbb{G}$, $\mathbb{G}'$ for which there exists a morphism $(F,\mu)$ which is invertible at least up to a 2-isomorphism. We are only interested in 2-groups up to equivalence.

\subsection{Classification up to equivalence}
\label{classificacio}

A basic result on 2-groups is Sinh's theorem \cite{hxSi75}. It says that {\it any} 2-group is equivalent to some sort of ``twisted'' version of a split 2-group $A[1]\rtimes G[0]$ for some $G$-module $A$. Because of its importance for what follows, we recall in this section the precise statement and how an equivalent ``twisted'' split 2-group is obtained from an arbitrary 2-group $\mathbb{G}$.

Let $\pi_0(\mathbb{G})$ be the set of isomorphism classes of objects in $\mathbb{G}$. It is a group with the product induced by the tensor product, i.e.
$$
[x]\cdot[y]=[x\otimes y].
$$
Let $\pi_1(\mathbb{G})$ be the group ${\rm Aut}_\Gg(I)$ of automorphisms of the unit object of $\mathbb{G}$. It is an abelian group with the composition (see \cite{SR72}), and it has a canonical $\pi_0(\mathbb{G})$-module structure given by
\begin{equation} \label{accio_pi0}
[x]\lhd u:=\gamma_x^{-1}(\delta_x(u)),
\end{equation}
for any $u\in\pi_1(\mathbb{G})$, any $[x]\in\pi_0(\mathbb{G})$ and any representative $x$ of $[x]$. Here $\delta_x,\gamma_x:\pi_1(\mathbb{G})\longrightarrow\mathrm{Aut}(x)$ stand for the canonical maps given by
\begin{align}
\delta_x(u)&=r_x\circ(u\otimes id_x)\circ r_x^{-1}, \label{delta_x}\\ \gamma_x(u)&=l_x\circ(id_x\otimes u)\circ l_x^{-1}, \label{gamma_x}
\end{align}
which are shown to be group isomorphisms. Hence it makes sense to consider the corresponding split 2-group $\pi_1(\mathbb{G})[1]\rtimes\pi_0(\mathbb{G})[0]$.

We are interested in this 2-group but equipped with a non-trivial associator
$$
a_{g,g',g''}:gg'g''\To gg'g''
$$
given by
\begin{equation} \label{associador_G}
a_{g,g',g''}=({\sf z}(g,g',g''),gg'g'')
\end{equation}
for some 3-cocycle ${\sf z}\in Z^3(\pi_0(\mathbb{G}),\pi_1(\mathbb{G}))$. This 3-cocycle is built from the associator of $\mathbb{G}$ as follows. Choose for each $g\in\pi_0(\mathbb{G})$ a representative $x_g\in g$, with $x_1=I$, and for each $y\in g$ an isomorphism $\iota_y:y\To x_g$, with $\iota_{x_g}=id_{x_g}$, and let
$$
{\sf z}:\pi_0(\mathbb{G})\times\pi_0(\mathbb{G})\times\pi_0(\mathbb{G})\To\pi_1(\mathbb{G})
$$
be the map uniquely defined by the commutativity of the diagrams
\begin{equation} \label{alpha}
\xymatrix{
x_{gg'g''}\ar[d]_{\iota^{-1}_{x_{gg'}\otimes x_{g''}}}\ar[rr]^{\gamma_{x_{gg'g''}}({\sf z}(g,g',g''))} && x_{gg'g''} \\ x_{gg'}\otimes x_{g''}\ar[d]_{\iota^{-1}_{x_g\otimes x_{g'}}\otimes id_{x_{g''}}} && x_g\otimes x_{g'g''}\ar[u]_{\iota_{x_g\otimes x_{g'g''}}} \\ (x_g\otimes x_{g'})\otimes x_{g''}\ar[rr]_{a_{x_g,x_{g'},x_{g''}}} & & x_g\otimes(x_{g'}\otimes x_{g''})\ar[u]_{id_{x_g}\otimes\iota_{x_{g'}\otimes x_{g''}}}
}
\end{equation}
for all $g,g',g''\in\pi_0(\G)$. It follows from the coherence conditions on the associator $a$ of $\mathbb{G}$ that ${\sf z}$ indeed is a 3-cocycle of $\pi_0(\mathbb{G})$ with values in $\pi_1(\mathbb{G})$ and that (\ref{associador_G}) indeed defines an associator for $\pi_1(\mathbb{G})[1]\rtimes\pi_0(\mathbb{G})[0]$. We shall denote the 2-group so defined by $\pi_1(\mathbb{G})[1]\rtimes_{{\sf z}}\pi_0(\mathbb{G})[0]$.

The 3-cocycle {\sf z} and hence, also the 2-group $\pi_1(\mathbb{G})[1]\rtimes_{\sf z}\pi_0(\mathbb{G})[0]$ obviously depend on the choices of representative objects $x_g$ and isomorphisms $\iota_y$. However, different choices lead to cohomologous 3-cocycles ${\sf z}_1$ and ${\sf z}_2$, and the 2-groups $\pi_1(\mathbb{G})[1]\rtimes_{{\sf z}_1}\pi_0(\mathbb{G})[0]$ and $\pi_1(\mathbb{G})[1]\rtimes_{{\sf z}_2}\pi_0(\mathbb{G})[0]$ are equivalent. In fact, Sinh's theorem says that we have
\begin{equation} \label{equivalencia}
\pi_1(\mathbb{G})[1]\rtimes_{\sf z}\pi_0(\mathbb{G})[0]\simeq\G
\end{equation}
for any of the above 3-cocycles {\sf z}. An explicit equivalence is given by the functor
$$
F:\Gg_{\pi_1(\G),\pi_0(\G)}\To\Gg
$$
defined on objects $g\in\pi_0(\G)$ and morphisms $(u,g)\in\pi_1(\G)\times\pi_0(\G)$ by
\begin{equation} \label{functor_F}
F(g)=x_g,\qquad F(u,g)=\gamma_{x_g}(u)
\end{equation}
and with the monoidal structure $\mu$ given by
\begin{equation} \label{mu}
\mu_{g,g'}=\iota^{-1}_{x_g\otimes x_{g'}}:x_{gg'}\To x_g\otimes x_{g'}
\end{equation}
(see Appendix).

It follows that the equivalence class of $\G$ is completely given by the groups $\pi_0(\mathbb{G})$ and $\pi_1(\mathbb{G})$, called {\it homotopy groups} of $\mathbb{G}$, and the cohomology class $\alpha=[{\sf z}]\in H^3(\pi_0(\G),\pi_1(\G))$, called its {\it Postnikov invariant}. Any 3-cocycle ${\sf z}\in\alpha$ will be called a {\it classifying 3-cocycle} of $\mathbb{G}$.

Notice that, because of the presence of the isomorphisms $\iota_y$ in (\ref{alpha}), the Postnikov invariant of $\mathbb{G}$ may be nontrivial even when the associator of $\mathbb{G}$ is the identity. In particular, it can be nonzero even for strict 2-groups. This can only happen when the strict 2-group is non skeletal. In fact, split 2-groups are those which are equivalent to strict skeletal ones, and they are characterized by the fact that $\alpha=0$.

\section{The 2-category of chain complexes}

From now on, $\kb$ denotes an arbitrary commutative ring with unit and $\Aa$ an arbitrary $\kb$-linear abelian category (i.e. an abelian category whose hom-sets are $\kb$-modules and such that the composition law is $\kb$-bilinear).

\subsection{(Split) chain complexes}

Recall that a chain complex in $\Aa$ is a sequence of objects $A_k$ of $\Aa$, labelled by $k\in\mathbb{Z}$, together with morphisms 
$$
d_k:A_{k+1}\To A_k,\qquad k\in\mathbb{Z}
$$
such that $d_{k-1}\circ d_k=0$ for all $k\in\mathbb{Z}$. Such a complex will be denoted $A_{\bullet}$. For any chain complex $A_{\bullet}$ and any $k\in\mathbb{Z}$ one defines
$$
Z_k(A_{\bullet}):=\mathrm{Ker}\left(A_k\stackrel{d_{k-1}}{\longrightarrow}A_{k-1}\right),
$$
$$
B_k(A_{\bullet}):=\mathrm{Ker}\left(\mathrm{Coker}\left(A_{k+1}\stackrel{d_{k}}{\longrightarrow}A_k\right)\right).
$$
$Z_k(A_{\bullet})$ is called the object of $k$-{\it cycles} of $A_{\bullet}$ and $B_k(A_{\bullet})$ the object of $k$-{\it boundaries} of $A_{\bullet}$. Both are subobjects of $A_k$. In fact, the above condition $d^2=0$ implies that $B_k(A_{\bullet})$ is a subobject of $Z_k(A_{\bullet})$, and one defines the $k$-{\it homology} object of $A_{\bullet}$ by
$$
H_k(A_{\bullet}):=\mathrm{Coker}\left(B_k(A_{\bullet})\hookrightarrow Z_k(A_{\bullet})\right).
$$
$A_{\bullet}$ is called {\it acyclic} or an {\it exact sequence} when $H_k(A_{\bullet})=0$ for all $k\in\mathbb{Z}$.

\begin{ex} \label{complex_escindit}{\rm
For any sequences $\{X_k\}_{k\in\mathbb{Z}}$ and $\{Y_{k}\}_{k\in\mathbb{Z}}$ of objects in $\Aa$, a chain complex can be defined as follows. Take
$$
A_k:=X_k\oplus Y_k\oplus X_{k-1},\qquad k\in\mathbb{Z}
$$
and let $d_k:A_{k+1}\To A_k$ be the morphism given by the composition
$$
X_{k+1}\oplus Y_{k+1}\oplus X_k\stackrel{\pi}{\longrightarrow} X_k\stackrel{\iota}{\longrightarrow} X_k\oplus Y_k\oplus X_{k-1},
$$
where $\pi$ and $\iota$ stand for the projection and injection associated to the corresponding biproduct. In this case, the condition $d^2=0$ follows because $\pi_j\circ\iota_i=0$ when $i\neq j$, where $\pi_j$ and $\iota_i$ refer here to the same biproduct. We clearly have
\begin{align*}
B_k(A_{\bullet})&\cong X_k, \\ Z_k(A_{\bullet})&\cong X_k\oplus Y_k, \\ H_k(A_{\bullet})&\cong Y_k
\end{align*}
for any $k\in\mathbb{Z}$. }
\end{ex}
We shall be mostly concerned with chain complexes isomorphic to the previous one for some sequences $\{X_k\}_{k\in\mathbb{Z}}$ and $\{Y_{k}\}_{k\in\mathbb{Z}}$. They are called {\it split} chain complexes and are characterized by the existence of the so called {\it splitting maps}. More precisely, we have the following.

\begin{prop}
A chain complex $A_{\bullet}$ in a $\kb$-linear abelian category $\Aa$ is split if and only if there exists morphisms (called {\em splitting maps}) $s_k:A_k\To A_{k+1}$, $k\in\mathbb{Z}$, such that $d_k\circ s_k\circ d_k=d_k $ for all $k\in\mathbb{Z}$.
\end{prop}
For short exact sequences $0\To A'\stackrel{f}{\To} A\stackrel{g}{\To} A''\To 0$ this criterion reduces to the existence of a section of $g$ (a map $s:A''\To A$ such that $g\circ s=id_{A''}$) or equivalently, a retraction of $f$ (a map $r:A\To A'$ such that $r\circ f=id_{A'}$), and in this case the original sequence is isomorphic to the split sequence $0\To A'\To A'\oplus A''\To A''\To 0$.

As it is well known, for some categories $\Aa$ all chain complexes are split. Examples are given by the category $\Vv ect_{\mathbb{F}}$ of vector spaces over an arbitrary field $\mathbb{F}$ and, more generally, any {\it semisimple} $\kb$-linear abelian category $\Aa$ (a $\kb$-linear abelian category such that all short exact sequences split). But this is not true in general. For instance, the chain complex (actually, an exact sequence) of abelian groups ($\kb=\mathbb{Z}$)
\begin{equation} \label{complex_grups_abelians}
\xymatrix{ \cdots\ar[r]^{\cdot 2} & \mathbb{Z}_4\ar[r]^{\cdot 2} & \mathbb{Z}_4\ar[r]^{\cdot 2} & \mathbb{Z}_4\ar[r]^{\cdot 2} & \cdots }
\end{equation}
does not split (for any morphism $s:\mathbb{Z}_4\To\mathbb{Z}_4$ we have $d\circ s\circ d=0$).

\subsection{Elementary description of the 2-category of chain complexes}

As pointed out by Gabriel and Zisman \cite{GZ67}, the chain complexes in an arbitrary $\kb$-linear~\footnote{In fact, Gabriel and Zisman consider the case $\kb=\mathbb{Z}$, but the generalization to arbitrary $\kb$ is straightforward.} abelian category $\Aa$ are the objects of a strict 2-category $\mathbf{Ch}(\Aa)$. In elementary terms, this is the 2-category described as follows. For a more conceptual description, which also paves the way toward the definition of an $\infty$-category of chain complexes, the reader is referred to \cite{GZ67}.

\paragraph{\em 1-morphisms.}
For arbitrary chain complexes $A_{\bullet}$ and $B_{\bullet}$, a 1-morphism from $A_{\bullet}$ to $B_{\bullet}$ is a morphism of chain complexes, i.e. a sequence of morphisms $f=\{f_k\}_{k\in\mathbb{Z}}$ in $\Aa$, with 
$$
f_k:A_k\To B_k,\qquad k\in\mathbb{Z},
$$
such that the diagram
$$
\xymatrix{
\cdots\ar[r] & A_{k+1}\ar[d]_{f_{k+1}}\ar[r]^{d_k} & A_k\ar[d]_{f_k}\ar[r]^{d_{k-1}} & A_{k-1}\ar[d]^{f_{k-1}}\ar[r] & \cdots \\ \cdots\ar[r] & B_{k+1}\ar[r]_{d_k} & B_k\ar[r]_{d_{k-1}} & B_{k-1}\ar[r] & \cdots }
$$
commutes. Among all 1-morphisms, the so called {\it null homotopic} ones play an special role. They are the 1-morphisms $f$ of the form
$$
f_k=d_{k}^B\circ h_k+h_{k-1}\circ d_{k-1}^A,\qquad k\in\mathbb{Z}
$$
for some family $h=\{h_k:A_k\To B_{k+1},k\in\mathbb{Z}\}$ of morphisms in $\Aa$. Such a family is called a chain contraction of $f$. In general, any family $h$ of morphisms as above will be called a {\it degree 1 homotopy} between $A_{\bullet}$ and $B_{\bullet}$, and two 1-morphisms $f,f'$ are called {\it homotopic} when their diference $f-f'$ is null homotopic.

\paragraph{\em 2-morphisms.}
For any parallel 1-morphisms $f,g:A_{\bullet}\longrightarrow B_{\bullet}$, a 2-morphism between them is the homotopy class of a chain homotopy between $f$ and $\tilde{f}$, i.e. of a degree 1 homotopy $h=\{h_k,k\in\mathbb{Z}\}$ between $A_{\bullet}$ and $B_{\bullet}$ such that
\begin{equation} \label{codomini_2-morfisme}
\tilde{f}_k=f_k+d_k^B\circ h_k+h_{k-1}\circ d_{k-1}^A,\qquad k\in\mathbb{Z}.
\end{equation}
Recall that two such degree 1 homotopies $h$ and $h'$ are said to be homotopic if and only if there exists a {\it homotopy of degree 2} between $A_{\bullet}$ and $B_{\bullet}$ relating them, i.e. a sequence of morphisms $h^{(2)}=\{h^{(2)}_k,k\in\mathbb{Z}\}$ in $\Aa$, with $h^{(2)}_k:A_k\To B_{k+2}$, such that
$$
h'_k=h_k+d_{k+1}^B\circ h^{(2)}_k-h^{(2)}_{k-1}\circ d_{k-1}^A,\qquad k\in\mathbb{Z}.
$$
Notice that the same sequence of morphisms $h_k:A_k\To B_{k+1}$ will be a 2-morphism between many parallel 1-morphisms, because the domain can be chosen arbitrarily. In this sense, 2-morphisms in $\mathbf{Ch}(\Aa)$ actually correspond to pairs $(f,[h])$, with $f$ any 1-morphism and $h$ any degree 1 homotopy. Then $f$ is the domain, while the codomain is given by (\ref{codomini_2-morfisme}).

\paragraph{\em Vertical composition of  2-morphisms.}

Given 2-morphisms
$$
\xymatrix{
A_{\bullet}\ar[rrrr]^f\ar@2{.}[d] && \ar@2{->}[d]^{(f,[h])} & & B_{\bullet}\ar@2{.}[d] \\ A_{\bullet}\ar[rrrr]^{\tilde{f}\hspace{0.6 truecm} }\ar@2{.}[d] && \ar@2{->}[d]^{(\tilde{f},[\tilde{h}])} && B_{\bullet}\ar@2{.}[d] \\ A_{\bullet}\ar[rrrr]_{\tilde{\tilde{f}}} &&& & B_{\bullet}
}
$$
the vertical composite is the 2-morphism
\begin{equation} \label{composicio_vertical_2-morfismes_1}
(\tilde{f},[\tilde{h}])\cdot(f,[h]):=(f,[h+\tilde{h}]),
\end{equation}
where $h+\tilde{h}$ is the degree 1 homotopy with components
\begin{equation} \label{composicio_vertical_2-morfismes_2}
(h+\tilde{h})_k=h_k+\tilde{h}_k.
\end{equation}
Observe that any 2-morphism is actually a 2-isomorphism, with inverse the 2-morphism described by the degree 1 homotopy having the same components but with the opposite sign.

\paragraph{\em Identity 2-morphisms.}

The identity 2-morphism of $f:A_{\bullet}\To B_{\bullet}$ is the pair $(f,[0])$, where $0$ stands for the homotopy of degree 1 having all components equal to the corresponding zero morphism in $Hom_{\Aa}(A_k,B_{k+1})$.

\paragraph{\em Composition between 1-morphisms.}

It is given by the usual composition of morphisms of complexes. Thus, for any 1-morphisms $f:A_{\bullet}\To B_{\bullet}$ and $g:B_{\bullet}\To C_{\bullet}$, the composite $g\circ f:A_{\bullet}\To C_{\bullet}$ is the 1-morphism with components
$$
(g\circ f)_k=g_k\circ f_k,\qquad k\in\mathbb{Z}.
$$

\paragraph{\em Horizontal composition of 2-morphisms.}

Suppose we are given 2-morphisms
$$
\xymatrix{
A_{\bullet}\ar[rrrr]^f\ar@2{.}[d] && \ar@2{->}[d]^{(f,[h])} & & B_{\bullet}\ar@2{.}[d]\ar[rrrr]^g\ar@2{.}[d] && \ar@2{->}[d]^{(g,[h'])} & & C_{\bullet}\ar@2{.}[d] \\ A_{\bullet}\ar[rrrr]_{\tilde{f}} &&&& B_{\bullet}\ar[rrrr]_{\tilde{g}} &&&& C_{\bullet}
}
$$
so that
\begin{align*}
\tilde{f}_k&=f_k+d_B\circ h_k+h_{k-1}\circ d_A \\ \tilde{g}_k&=g_k+d_C\circ h'_k+h'_{k-1}\circ d_B
\end{align*}
for all $k\in\mathbb{Z}$. Then the composite is given by
\begin{equation} \label{composicio_horitzontal_2-morfismes_1}
(g,[h'])\circ(f,[h]):=(g\circ f,[\hat{h}])
\end{equation}
where $\hat{h}$ is the homotopy of degree 1 with components
\begin{equation} \label{composicio_horitzontal_2-morfismes_2}
\hat{h}_k:=h'_k\circ f_k+h'_k\circ d_B\circ h_k+h'_k\circ h_{k-1}\circ d_A+g_{k+1}\circ h_k.
\end{equation}

\begin{rem}{\rm 
In fact, there are two ways of thinking of the composition $\tilde{g}_k\circ\tilde{f}_k$ as the sum of $g_k\circ f_k$ plus a homotopically trivial term. This leads to two different expressions for the homotopy of first degree $\hat{h}$. Indeed, the above composition gives
\begin{align*}
\tilde{g}_k\circ\tilde{f}_k&=g_k\circ f_k+g_k\circ d_B\circ h_k+g_k\circ h_{k-1}\circ d_A \\ &\ \ \ \ \ \ +d_C\circ h'_k\circ f_k+d_C\circ h'_k\circ d_B\circ h_k+d_C\circ h'_k\circ h_{k-1}\circ d_A \\ &\ \ \ \ \ \ +h'_{k-1}\circ d_B\circ f_k+h'_{k-1}\circ d_B\circ d_B\circ h_k+h'_{k-1}\circ d_B\circ h_{k-1}\circ d_A \\ &=g_k\circ f_k+d_C\circ g_{k+1}\circ h_k+g_k\circ h_{k-1}\circ d_A \\ &\ \ \ \ \ \ +d_C\circ h'_k\circ f_k+d_C\circ h'_k\circ d_B\circ h_k+d_C\circ h'_k\circ h_{k-1}\circ d_A \\ &\ \ \ \ \ \ +h'_{k-1}\circ f_{k+1}\circ d_A+h'_{k-1}\circ d_B\circ h_{k-1}\circ d_A.
\end{align*}
The dilemma arises with the term $d_C\circ h'_k\circ h_{k-1}\circ d_A$, which can be groupped either with the terms of the form $d_C\circ -$ or with those of the form $-\circ d_A$. The expression (\ref{composicio_horitzontal_2-morfismes_2}) for $\hat{h}$ corresponds to the first option, while for the second option it is given by
$$
\hat{h}_k=g_{k+1}\circ h_k+h'_k\circ f_k+h'_k\circ d_B\circ h_k+d_C\circ h'_{k+1}\circ h_k.
$$
Both expressions are, however, homotopic and hence, lead to the same composite 2-morphism.
}
\end{rem}

\section{The 2-group of symmetries of a split chain complex}

\subsection{General definition and examples}
\label{definicio-exemples}

Let $A_{\bullet}$ be any chain complex in $\Aa$, split or not. According to the previous description of the 2-category $\mathbf{Ch}(\Aa)$, its 2-group of symmetries $\mathbb{E}quiv(A_{\bullet})$ is given as follows:
\begin{itemize}
\item objects: the self-equivalences of $A_{\bullet}$, i.e. 1-endomorphisms $f:A_{\bullet}\To A_{\bullet}$ for which there exists another 1-endomorphism $f^*:A_{\bullet}\To A_{\bullet}$ and degree 1 homotopies $h^{(1)}$ and $\tilde{h}^{(1)}$ from $A_{\bullet}$ to itself such that
\begin{align*}
f^*_k\circ f_k&=id_{A_k}+d_k\circ h^{(1)}_k+h^{(1)}_{k-1}\circ d_{k-1} \\ f_k\circ f^*_k&=id_{A_k}+d_k\circ\tilde{h}^{(1)}_k+\tilde{h}^{(1)}_{k-1}\circ d_{k-1}
\end{align*}
for all $k\in\mathbb{Z}$;
\item morphisms: the 2-morphisms between self-equivalences; hence, for any self-equivalences $f,f'$ of $A_{\bullet}$ the hom-set $\mathrm{Hom}_{\Ee quiv(A_{\bullet})}(f,f')$ is the set of pairs $(f,[h^{(1)}])$ with $h^{(1)}:A_{\bullet}\To A_{\bullet}$ any degree 1 homotopy such that
$$
f'_k=f_k+d_k\circ h^{(1)}_k+h^{(1)}_{k-1}\circ d_{k-1},\qquad \forall\ k\in\mathbb{Z}
$$
and $[h^{(1)}]$ the corresponding homotopy class (in particular, all morphisms are invertible because all 2-morphisms in $\mathbf{Ch}(\Aa)$ are invertible);
\item composition law: it is given by the vertical composition of 2-morphisms in $\mathbf{Ch}(\Aa)$, i.e.
\begin{equation} \label{composicio}
(f',[\overline{h}^{(1)}])\circ(f,[h^{(1)}])=(f,[\overline{h}^{(1)}+h^{(1)}]),
\end{equation}
where $f'$ denotes the codomain of $(f,[h^{(1)}])$;
\item identity morphisms: for any object $f$, its identity morphism is $id_f=(f,[0])$;

\item tensor product: it is given on objects by the composition of 1-morphisms, and on morphism by the horizontal composition of 2-morphisms; thus for any self-equivalences $f,\tilde{f}$ of $A_{\bullet}$ we have
\begin{align}
\tilde{f}\otimes f&=\tilde{f}\circ f, \label{prod_tensorial_1} \\ (\tilde{f},[\tilde{h}^{(1)}])\otimes(f,[h^{(1)}])&=(\tilde{f}\circ f,[\hat{h}^{(1)}]), \label{prod_tensorial_2}
\end{align}
with $\hat{h}^{(1)}$ the degree 1 homotopy with components given by
\begin{equation} \label{hat_h}
\hat{h}^{(1)}_k=\tilde{h}^{(1)}_k\circ f_k+\tilde{h}^{(1)}_k\circ d_B\circ h^{(1)}_k+\tilde{h}^{(1)}_k\circ h^{(1)}_{k-1}\circ d_A+\tilde{f}_{k+1}\circ h^{(1)}_k;
\end{equation}
\item unit object: $id_{A_{\bullet}}$.
\end{itemize}
Observe that the underlying groupoid $\mathcal{E}quiv(A_{\bullet})$ of $\mathbb{E}quiv(A_{\bullet})$ is strict because the 2-category of chain complexes is strict. However, the 2-group itself will be non-strict in general, because there may exists self-equivalences of $A_{\bullet}$ which are not automorphisms. We are only interested in the equivalence class of $\mathbb{E}quiv(A_{\bullet})$. This is completely determined by the corresponding homotopy groups, which we shall denote $\pi_0(A_\bullet)$ and $\pi_1(A_\bullet)$, and the Postnikov invariant.

\begin{ex} \label{ex1}{\rm
For any $k\in\mathbb{Z}$, $k\neq 0$, let $A_{\bullet}$ be the {\it non-split} chain complex of abelian groups
$$
\xymatrix{
0\ar[r] & \mathbb{Z}\ar[r]^{\cdot 2k} & \mathbb{Z}\ar[r] & \mathbb{Z}_2\ar[r] & 0\ , }
$$
which we shall assume concentrated in degrees 0, 1 and 2 (no splitting maps exist because the unique morphism of groups from $\mathbb{Z}_2$ to $\mathbb{Z}$ is the zero one). Then $\pi_1(A_\bullet)$ is trivial, $\pi_0(A_\bullet)$ is isomorphic to $\mathbb{Z}^*_{2k}$ and the Postnikov invariant is zero, so that
$$
\mathbb{E}quiv(A_{\bullet})\simeq \mathbb{Z}^*_{2k}[0].
$$
Here $\mathbb{Z}^*_{2k}$ denotes the (multiplicative) group of units of the ring $\mathbb{Z}_{2k}$. To see this, let us first observe that an arbitrary endomorphism of $A_{\bullet}$ is either of the form
$$
\xymatrix{
0\ar[r] & \mathbb{Z}\ar[r]^{\cdot 2k}\ar[d]_{\cdot n} & \mathbb{Z}\ar@{.>}[dl]_h\ar[r]\ar[d]^{\cdot n} & \mathbb{Z}_2\ar[r]\ar[d]^{id} & 0 \\ 0\ar[r] & \mathbb{Z}\ar[r]^{\cdot 2k} & \mathbb{Z}\ar[r] & \mathbb{Z}_2\ar[r] & 0 }
$$
for odd $n$, or of the form
$$
\xymatrix{
0\ar[r] & \mathbb{Z}\ar[r]^{\cdot 2k}\ar[d]_{\cdot n} & \mathbb{Z}\ar@{.>}[dl]_h\ar[r]\ar[d]^{\cdot n} & \mathbb{Z}_2\ar[r]\ar[d]^{0} & 0 \\ 0\ar[r] & \mathbb{Z}\ar[r]^{\cdot 2k} & \mathbb{Z}\ar[r] & \mathbb{Z}_2\ar[r] & 0 }
$$
if $n$ is even. Let us denote this endomorphism by $f^{(n)}$ for any $n\geq\mathbb{Z}$. In particular, we have $f^{(1)}=id_{A_{\bullet}}$. Now, for any $n,m\in\mathbb{Z}$ we have
\begin{equation} \label{homotops}
[f^{(n)}]=[f^{(m)}]\ \ \Leftrightarrow\ \ n\equiv m \ (\mathrm{mod}\ 2k).
\end{equation}
Indeed, $f^{(n)}$ and $f^{(m)}$ are homotopic if and only if there exists a morphism of groups $h:\mathbb{Z}\To\mathbb{Z}$ (as shown in the above diagrams), hence a map of the form $h=\phi^{(l)}$ for some $l\in\mathbb{Z}$, with $\phi^{(l)}$ the endomorphism of $\mathbb{Z}$ uniquely defined by $\phi^{(l)}(1)=l$, such that
\begin{align*}
f^{(n)}_0&=f^{(m)}_0 \\ \phi^{(n)}&=\phi^{(m)}+\phi^{(2k)}\circ\phi^{(l)}.
\end{align*} 
But this is equivalent to the equality
$$
n=m+2kl
$$
in $\mathbb{Z}$, from which (\ref{homotops}) follows (the first condition says that $n,m$ must be of the same parity and hence, it is subsumed by this condition). In particular, we have
$$
f^{(n)}\simeq id_{A_{\bullet}}\ \Leftrightarrow\ n\equiv 1\ \mathrm{mod}~2k.
$$
Since $f^{(n')}\circ f^{(n)}=f^{(n'n)}$, it follows that $f^{(n)}$ is a self-equivalence of $A_{\bullet}$ if and only if $n\in\mathbb{Z}^*_{2k}$. In summary, we have
$$
\mathrm{Equiv}(A_{\bullet})=\{f^{(n)},\ n\in\mathbb{Z}\ \mbox{such that}\ hcf(n,2k)=1\}
$$
(this indeed is a submonoid of the multiplicative monoid $\mathrm{End}(A_{\bullet})$ of all endomorphisms), and
$$
\pi_0(A_{\bullet})=\{[f^{(n)}],\ n\in\{0,1,\ldots,2k-1\}\ \mbox{such that}\ [n]\in\mathbb{Z}^*_{2k}\}\ \cong\ \mathbb{Z}^*_{2k},
$$
the isomorphism being as (multiplicative) groups. Finally, $\pi_1(A_{\bullet})$ is trivial because there is a unique degree 1 homotopy from $id_{A_{\bullet}}$ to itself, namely the zero one. Indeed, any such homotopy is completely given by an endomorphism $\phi^{(l)}:\mathbb{Z}\To\mathbb{Z}$ such that $id_{\mathbb{Z}}=id_{\mathbb{Z}}+\phi^{(2k)}\circ\phi^{(l)}$ and hence, such that $1=1+2kl$ in $\mathbb{Z}$, from which it follows $l=0$.
}
\end{ex}

\begin{ex} \label{ex2}{\rm
Let $A_{\bullet}$ be the {\it non-split exact sequence} of abelian groups (\ref{complex_grups_abelians}). Then $\pi_0(A_\bullet)$ is trivial, $\pi_1(A_\bullet)$ is isomorphic to $\mathbb{Z}_{2}$ and the Postnikov invariant is zero, so that
$$
\mathbb{E}quiv(A_{\bullet})\simeq \mathbb{Z}_2[1].
$$
To prove this, let us first observe that any endomorphism $f=\{f_k\}_{k\in\mathbb{Z}}$ of (\ref{complex_grups_abelians}) has a well defined parity in the following sense. For any $l\in\mathbb{Z}_4$, let us denote by $\psi^{(l)}$ the unique endomorphism of $\mathbb{Z}_4$ such that $\psi^{(l)}(1)=l$. In particular, we have $\psi^{(0)}=0$, $\psi^{(1)}=id_{\mathbb{Z}_4}$ and $d_k=\psi^{(2)}$ for all $k\in\mathbb{Z}$. Notice also that
\begin{align*}
\psi^{(l)}\circ\psi^{(l')}&=\psi^{(ll')} \\ \psi^{(l)}+\psi^{(l')}&=\psi^{(l+l')}.
\end{align*}
Then if for some $k\in\mathbb{Z}$ we have $f_k=\psi^{(l)}$ it follows from the morphism condition that both $f_{k-1}$ and $f_{k+1}$ are necessarily equal to $\psi^{(l)}$ or $\psi^{(l+2)}$. Hence all components of $f$ are of the form $\psi^{(l)}$ with either all $l$ odd or all $l$ even. In case $f$ is odd, which components $f_k$ are equal to $\psi^{(1)}$ and which are equal to $\psi^{(3)}$ can be chosen arbitrarily, and similarly when $f$ is even. In fact, we have a morphism of monoids
$$
\epsilon:\mathrm{End}(A_{\bullet})\To\mathbb{Z}_2
$$
defined by
$$
\epsilon(f)=\left\{\begin{array}{ll} 0,&\mbox{if $f$ is even} \\ 1, & \mbox{if $f$ is odd} \end{array}\right.
$$
(here we think of $\mathbb{Z}_2$ as a {\it multiplicative} monoid). Next observation is that
$$
f\simeq f'\ \Leftrightarrow\ \epsilon(f)=\epsilon(f').
$$
Indeed, let $f=\{\psi^{(l_k)}\}_{k\in\mathbb{Z}}$ and  $f'=\{\psi^{(l'_k)}\}_{k\in\mathbb{Z}}$. Then it is immediate to check that $f$ is homotopic to $f'$ if and only if there exists a sequence $\{\tilde{l}_k\}_{k\in\mathbb{Z}}$ of elements of $\mathbb{Z}_4$ such that
$$
l'_k=l_k+2(\tilde{l}_k+\tilde{l}_{k-1}),\qquad \forall\ k\in\mathbb{Z}
$$
(the sequence $\{\tilde{l}_k\}$ gives the components $h_k=\psi^{(\tilde{l}_k)}$ of a degree 1 homotopy $h$ between $f$ and $f'$). Now if such a sequence exists, both $l_k$ and $l'_k$ are clearly of the same parity (hence $\epsilon(f)=\epsilon(f')$) and conversely, if $l_k,l'_k$ are of the same parity, a sequence $\{\tilde{l}_k\}$ as required can be built (non uniquely) starting with $\tilde{l}_0=0$, for example, and applying the above recurrence to compute $\tilde{l}_k$ for all $k>0$ and for all $k<0$ separately. It follows that $f$ is a self-equivalence of $A_{\bullet}$ if and only if there exists $f'$ such that
$$
\epsilon(f')\epsilon(f)=1
$$
and therefore, if and only if $\epsilon(f)=1$. We conclude that
$$
\mathrm{Equiv}(A_{\bullet})=\{f\in\mathrm{End}(A_{\bullet})\ |\ \epsilon(f)=1\}
$$
and
$$
\pi_0(A_{\bullet})=\{[id_{A_{\bullet}}]\}\cong 1.
$$
Finally, it is easy to check that any degree 1 homotopy $h$ from $id_{A_{\bullet}}$ to itself is given by a sequence of endomorphisms $\psi^{(\tilde{l}_k)}$ of $\mathbb{Z}_4$ with all $\tilde{l}_k\in\mathbb{Z}_4$ also of the same parity. Thus if $H(id_{A_{\bullet}},id_{A_{\bullet}})$ denotes the set of all such homotopies, we have a morphism of groups 
$$
\varepsilon:H(id_{A_{\bullet}},id_{A_{\bullet}})\To\mathbb{Z}_2
$$
mapping any homotopy to 0, if it is even, or to 1, if it is odd (now $\mathbb{Z}_2$ is thought of as an {\it additive} group; cf. (\ref{composicio_vertical_2-morfismes_2})). Moreover, two such homotopies $h,h'$ turn out to be homotopic if and only if $\varepsilon(h)=\varepsilon(h')$, so that
$$
\pi_1(A_{\bullet})\cong\mathbb{Z}_2
$$
as claimed before.
}
\end{ex}

\subsection{Case of an arbitrary chain complex}
\label{cas_general}

Our goal is to compute the homotopy groups and Postnikov invariant of $\mathbb{E}quiv(A_{\bullet})$ when $A_{\bullet}$ is split. Now, for an arbitrary chain complex $A_{\bullet}$, not necessarily split, these invariants are related to the {\it complexes of homotopies} between certain chain complexes. Hence we begin by recalling the definition of these complexes of homotopies (see \cite{GZ67}) and explaining their relationship to the invariants that classify $\mathbb{E}quiv(A_{\bullet})$ in the generic case.

\paragraph{\em The complex of $\kb$-modules $Hom(A_{\bullet},B_{\bullet})$.}
\label{seccio_morfisme_canonic}
For any chain complexes $A_{\bullet}$ and $B_{\bullet}$ in $\Aa$, the {\it complex of homotopies} between them is the complex of $\kb$-modules $Hom(A_{\bullet},B_{\bullet})$ whose piece in degree $k$ is
$$
Hom(A_{\bullet},B_{\bullet})_k:=\prod_{k'\in\mathbb{Z}}\Aa(A_{k'},B_{k'+k}),\qquad k\in\mathbb{Z},
$$
and whose boundary operator $d_{k-1}:Hom(A_{\bullet},B_{\bullet})_k\longrightarrow Hom(A_{\bullet},B_{\bullet})_{k-1}$ is given by
\begin{equation} \label{diferencial_homs}
d(h^{(k)})_{k'}=d_{k'+k-1}^B\circ h^{(k)}_{k'}-(-1)^kh^{(k)}_{k'-1}\circ d_{k'-1}^A,\qquad k'\in\mathbb{Z},
\end{equation}
for any $h^{(k)}\in Hom(A_{\bullet},B_{\bullet})_k$. The elements of $Hom(A_{\bullet},B_{\bullet})_k$ will be called {\it degree $k$ homotopies} between $A_{\bullet}$ and $B_{\bullet}$ (or between the underlying $\mathbb{Z}$-graded objects). We already met homotopies of degree 1 and 2 when describing the 2-morphisms in the 2-category of chain complexes. As done before, we sometimes use a superscript to indicate the degree of a homotopy.

For our purposes we shall be concerned only with the degree zero component $Hom(A_{\bullet},B_{\bullet})_0$ of this complex~\footnote{This is because we are only interested in the {\it 2-category} of chain complexes and the {\it 2-groups} of symmetries of its objects. The whole complex will be necessary in order  to describe the $\infty$-{\it category} of chain complexes and the $\infty$-{\it groups} of symmetries of its objects.}, i.e. the $\kb$-module of all sequences $(h^{(0)}_{k'})_{k'\in\mathbb{Z}}$ of arbitrary morphisms $h^{(0)}_{k'}:A_{k'}\To B_{k'}$. It follows from (\ref{diferencial_homs}) that the associated object of 0-cycles $Z_0(Hom(A_{\bullet},B_{\bullet}))$ (resp. 0-boundaries $B_0(Hom(A_{\bullet},B_{\bullet}))$) is the sub-$\kb$-module of all (1-)morphisms of complexes (resp. null homotopic (1-)morphisms of complexes) between $A_{\bullet}$ and $B_{\bullet}$. Therefore the elements of $H_0(Hom(A_{\bullet},B_{\bullet}))$ are nothing but the homotopy classes of morphisms of complexes between $A_{\bullet}$ and $B_{\bullet}$.

Let us also recall that any morphism of chain complexes $f:A_{\bullet}\To B_{\bullet}$ induces a canonical $\kb$-linear map
$$
\Phi:H_0(Hom(A_{\bullet},B_{\bullet}))\longrightarrow \prod_{k\in\mathbb{Z}}Hom_{\Aa}(H_k(A_{\bullet}),H_k(B_{\bullet}))
$$
given by
\begin{equation} \label{morfisme_canonic}
\Phi([f]):=(H_k(f))_{k\in\mathbb{Z}},
\end{equation}
where $H_k(f):H_k(A_{\bullet})\To H_k(B_{\bullet})$ are the morphisms induced in homology. When $A_{\bullet}=B_{\bullet}$, $End(A_{\bullet})_0$ is a $\kb$-algebra with the product given by componentwise composition, and this structure is inherited by both $Z_0(End(A_{\bullet}))$ and $H_0(End(A_{\bullet}))$. Indeed, it is easy to check that $B_0(End(A_{\bullet}))$ is an ideal of $Z_0(End(A_{\bullet}))$. Moreover, by functoriality of the homology, $\Phi$ is in fact a morphism of $\kb$-algebras in this case.

In general, $\Phi$ is neither injective (non-homotopic morphisms may induce the same morphisms on homology) nor surjective (there may exist sequences of morphisms $\varphi_k:H_k(A_{\bullet})\To H_k(B_{\bullet})$, $k\in\mathbb{Z}$, which do not come from a morphism of complexes between $A_{\bullet}$ and $B_{\bullet}$). As shown below, however, it is an isomorphism of $\kb$-algebras when $A_{\bullet}=B_{\bullet}$ is a split chain complex.

\paragraph{\em The group  $\pi_0(A_{\bullet})$.}

$\pi_0(A_{\bullet})$ is the group of homotopy classes of self-equivalences of $A_{\bullet}$, with the product induced by the composition of self-equivalences. Now, if $[f]$ denotes the homotopy class of an endomorphism $f:A_{\bullet}\To A_{\bullet}$, $f$ is a self-equivalence if and only if there exists $f^*:A_{\bullet}\To A_{\bullet}$ such that
$$
[f^*][f]=[id_{A_{\bullet}}]=[f][f^*]
$$
in $H_0(End(A_{\bullet}))$. Therefore
\begin{equation} \label{pi_0_grup_unitats}
\pi_0(A_{\bullet})=U(H_0(End(A_{\bullet}))),
\end{equation}
where $U(H_0(End(A_{\bullet})))$ denotes the group of units of the $\kb$-algebra $H_0(End(A_{\bullet}))$.

\paragraph{\em The abelian group $\pi_1(A_{\bullet})$.}

Since all 2-morphisms in $\mathbf{Ch}(\Aa)$ are invertible, $\pi_1(A_{\bullet})$ is the abelian group of {\it all} 2-endomorphisms of $id_{A_{\bullet}}$. Now, according to (\ref{codomini_2-morfisme}), for any 2-morphism $(id_{A_{\bullet}},[h^{(1)}])$ of domain $id_{A_{\bullet}}$ its codomain is the morphism $f$ with components
$$
f_k=id_{A_k}+d_k\circ h^{(1)}_k+h^{(1)}_{k-1}\circ d_{k-1},\qquad k\in\mathbb{Z}.
$$
It follows that $(id_{A_{\bullet}},[h^{(1)}])$ is a 2-endomorphism of $id_{A_{\bullet}}$ if and only if the degree 1 homotopy $h^{(1)}$ is such that
\begin{equation} \label{automorfismes_id}
d_k\circ h^{(1)}_k+h^{(1)}_{k-1}\circ d_{k-1}=0
\end{equation}
for all $k\in\mathbb{Z}$ or, equivalently, if and only if $h^{(1)}$ is a morphism of chain complexes between $A_{\bullet}$ and its translation $A[1]_{\bullet}$, defined by
$$
A[1]_k=A_{k+1},\quad  d_k^{A[1]}=-d_{k+1}^A,\qquad k\in\mathbb{Z}.
$$
Moreover, degree 1 homotopies between $A_{\bullet}$ and $A[1]_{\bullet}$ are exactly the same as degree 2 homotopies from $A_{\bullet}$ into itself, and the relation ``being homotopic'' is exactly the same in both sets of homotopies, as the reader may easily check. Therefore taking the homotopy class of $h^{(1)}$ either as an element of
$$
Z_0(Hom(A_{\bullet},A[1]_{\bullet})\subset Hom(A_{\bullet},A[1]_{\bullet})_0
$$
or as an element of $End(A_{\bullet})_1$ gives exactly the same, and we have
\begin{equation} \label{pi_1_morfismes_A_A1}
\pi_1(A_{\bullet})=H_0(Hom(A_{\bullet},A[1]_{\bullet})).
\end{equation}
In fact, the equality is as abelian groups, because the ``sum'' of $\pi_1(A_{\bullet})$ is given by the vertical composition of 2-morphisms in $\mathbf{Ch}(\Aa)$, and this indeed corresponds to summing homotopies (cf. (\ref{composicio_vertical_2-morfismes_1})-(\ref{composicio_vertical_2-morfismes_2})).

\begin{rem} \label{remarca_infinit-grup}{\rm
In fact, the symmetries of an arbitrary chain complex $A_{\bullet}$ are expected to be the objects of an $\infty$-{\it group} (i.e. a one-object $\infty$-groupoid) whose homotopy groups will be given by
$$
\pi_n(A_{\bullet})=H_0(Hom(A_{\bullet},A[n]_{\bullet}))
$$
for all $n\geq 1$. 
}
\end{rem}

\paragraph{\em Structure of $\pi_0(A_{\bullet})$-module on $\pi_1(A_{\bullet})$.}

For a (possibly non-strict) 2-group $\mathbb{G}$ with underlying strict monoidal groupoid, the general expression (\ref{accio_pi0}) for the action of $\pi_0(\mathbb{G})$ on $\pi_1(\mathbb{G})$ reduces to
$$
[x]\lhd u=\epsilon\circ(id_{x^*}\otimes u\otimes id_x)\circ\epsilon^{-1},
$$
with $x^*$ any pseudoinverse of the chosen representative $x$ and $\epsilon:x^*\otimes x\stackrel{\cong}{\To} e$ any isomorphism (cf. Appendix). In particular, this is true for our 2-group $\mathbb{E}quiv(A_{\bullet})$, in which case $x$ is a self-equivalence $f$ of $A_{\bullet}$ and
$$
u=(id_{A_{\bullet}},[h^{(1)}]),
$$
with $h^{(1)}$ any degree 1 homotopy of $A_{\bullet}$ into itself satisfying (\ref{automorfismes_id}). Hence, if $f^*$ is any pseudoinverse of $f$ and $\overline{h}^{(1)}$ any degree 1 homotopy between $f^*\circ f$ and $id_{A_{\bullet}}$, we have
$$
[f]\lhd(id_{A_{\bullet}},[h^{(1)}])=(f^*\circ f,[\overline{h}^{(1)}])\circ((f^*,[0])\otimes(id_{A_{\bullet}},[h^{(1)}])\otimes(f,[0]))\circ(id_{A_{\bullet}},[-\overline{h}^{(1)}])=(id_{A_{\bullet}},[\hat{h}^{(1)}])
$$
with 
$$
\hat{h}^{(1)}_k=f^*_{k+1}\circ h^{(1)}_k\circ f_k,\qquad k\in\mathbb{Z}
$$
(see (\ref{composicio}), (\ref{prod_tensorial_2}) and (\ref{hat_h})). Observe that this is a homotopy which still satisfies (\ref{automorfismes_id}). Therefore, thinking of the elements of $\pi_1(A_{\bullet})$ as homotopy classes of morphisms $g:A_{\bullet}\To A[1]_{\bullet}$, the action of $\pi_0(A_{\bullet})$ on $\pi_1(A_{\bullet})$ is simply given by `conjugation'. More precisely, we have
\begin{equation} \label{accio_pi0_pi1}
[f]\lhd[g]=[f^*[1]\circ g\circ f],
\end{equation} 
where $f^*[1]:A[1]_{\bullet}\To A[1]_{\bullet}$ denotes the self-equivalence of $A[1]_{\bullet}$ induced by $f^*:A_{\bullet}\To A_{\bullet}$, with components $f^*[1]_k=f^*_{k+1}$ for all $k\in\mathbb{Z}$.

\paragraph{\em Postnikov invariant.}

Since the underlying monoidal groupoid of $\mathbb{E}quiv(A_{\bullet})$ is strict, its set of objects $Equiv(A_{\bullet})$ is a monoid with the product given by the tensor product and with $id_{A_{\bullet}}$ as unit. Moreover, the canonical projection
$$
\pi:Equiv(A_{\bullet})\To\pi_0(A_{\bullet})
$$
mapping any self-equivalence $f$ to its homotopy class $[f]$ is a morphism of monoids. The 2-group $\mathbb{E}quiv(A_{\bullet})$ will be split when this projection admits a section in the category of monoids. Indeed, let
$$
s:\pi_0(A_{\bullet})\To Equiv(A_{\bullet})
$$
be such a section. Take as representative of $[f]\in\pi_0(A_{\bullet})$ the self-equivalence $s[f]\in Equiv(A_{\bullet})$ and apply the algorithm described in \S~\ref{classificacio} to construct a classifying 3-cocycle $z$. We have that all $\iota$'s appearing in (\ref{alpha}) are identities because $s$ is a morphism of monoids, and hence
$$
z([f],[f'],[f''])=(id_{A_{\bullet}},[0])
$$
for all $[f],[f'],[f'']\in\pi_0(A_{\bullet})$.

As we shall see in the next paragraph, a section $s$ as before indeed exists when $A_{\bullet}$ is split. There also are non-split complexes, however, whose 2-group of symmetries are also split (for ex. the complexes in Examples~\ref{ex1} and \ref{ex2} above), even trivial up to equivalence (Example~\ref{ex1} with $k=2$).

\subsection{Case of a split chain complex}

It follows from the general 2-categorical yoga that equivalent objects in a 2-category have equivalent 2-groups of symmetries. Hence, for the sake of simplicity and without loss of generality, we shall assume from now on that $A_{\bullet}$ stands for the split chain complex of Example~\ref{complex_escindit}. To emphasize the fact that the objects $X_k$ and $Y_k$ respectively give the $k$-boundary and $k$-homology objects of $A_{\bullet}$, we shall denote them $B_k$ and $H_k$ respectively.

\begin{thm} \label{teorema_principal}
Let $A_{\bullet}$ be the split chain complex defined by the objects $\{B_k,k\in\mathbb{Z}\}$ and $\{H_k,k\in\mathbb{Z}\}$ of $\Aa$. Let $Aut_\Aa(H_{\bullet})$ be the group
$$
Aut_\Aa(H_{\bullet}):=\prod_{k\in\mathbb{Z}}Aut_{\Aa}(H_k),
$$
and let $Hom_{\Aa}(H_{\bullet},H_{\bullet+1})$ be the $\kb$-module
$$
Hom_{\Aa}(H_{\bullet},H_{\bullet+1}):=\prod_{k\in\mathbb{Z}}Hom_{\Aa}(H_k,H_{k+1})
$$
equipped with the $Aut_\Aa(H_{\bullet})$-module structure given by ``conjugation'', i.e.
$$
(\mbox{\boldmath$\psi$\unboldmath}\lhd\mbox{\boldmath$\xi$\unboldmath})_k=\psi^{-1}_{k+1}\circ\xi_k\circ\psi_k,\qquad k\in\mathbb{Z},
$$
for any $\mbox{\boldmath$\psi$\unboldmath}\in Aut_{\Aa}(H_\bullet)$ and any $\mbox{\boldmath$\xi$\unboldmath}\in Hom_{\Aa}(H_\bullet,H_{\bullet+1})$. Then we have an equivalence of 2-groups
\begin{equation} \label{equivalencia_final}
\mathbb{E}quiv(A_{\bullet})\simeq Hom_{\Aa}(H_\bullet,H_{\bullet+1})[1]\rtimes Aut_{\Aa}(H_\bullet)[0].
\end{equation}
In particular, the 2-group of symmetries of any split exact sequence is trivial (up to equivalence).
\end{thm}

\begin{ex}{\rm
Let $\mathbb{F}$ be a field and $d:V\To W$ an $\mathbb{F}$-linear map between arbitrary vector spaces over $\mathbb{F}$. We may think of $d$ as the chain complex concentrated in degrees 1 and 0 (as any chain complex in $\ev_\mathbb{F}$, it is split), and as such it has a split 2-group of symmetries given by
\begin{equation} \label{2-grup_lineal_general}
\mathbb{E}quiv\left(V\stackrel{d}{\To}W\right)\simeq Hom_{\kb}(\mathrm{Coker}~d,\mathrm{Ker}~d)[1]\rtimes(GL_{\kb}(\mathrm{Coker}~d)\times GL_{\kb}(\mathrm{Ker}~d))[0].
\end{equation}
In particular:
\begin{itemize}
\item If $d$ is monic, $\mathbb{E}quiv(d)$ is discrete with $GL_\kb(\mathrm{Coker}~d)$ as underlying group.
\item If $d$ is epi, $\mathbb{E}quiv(d)$ is discrete with $GL_\kb(\mathrm{Ker}~d)$ as underlying group.
\item If $d$ is an isomorphism, $\mathbb{E}quiv(d)$ is trivial (up to equivalence).
\end{itemize}
As discussed by Baez and Crans in \cite{BC03}, there is a sense in which $d$ can be considered a 2-vector space, i.e. a categorical analog of a vector space. Its 2-group of symmetries (\ref{2-grup_lineal_general}) then gives the corresponding {\it general linear 2-group}.}
\end{ex}

\noindent
The rest of this section is devoted to prove the above theorem. We shall first compute the homotopy groups, next we identify how the first acts on the second and finally, we show that the 2-group is split. At the end of the section, we show by explicit construction that there exists an equivalence (\ref{equivalencia_final}) which is given by a strict monoidal functor.

\paragraph{\em The group $\pi_0(A_{\bullet})$.}

According to (\ref{pi_0_grup_unitats}), $\pi_0(A_{\bullet})$ is the group of units of the $\kb$-algebra $H_0(End(A_{\bullet}))$. To identify this $\kb$-algebra, we first identify the $\kb$-algebra $End_{\Cc h(\Aa)}(A_{\bullet})$ of endomorphisms of $A_{\bullet}$ (i.e. the $\kb$-algebra $Z_0(End(A_{\bullet}))$ of 0-cycles of the complex $End(A_{\bullet})$), and then we take quotient by the homotopy relation.

By definition, the boundary map of $A_{\bullet}$
$$
d_k:B_{k+1}\oplus H_{k+1}\oplus B_k\longrightarrow B_k\oplus H_k\oplus B_{k-1},\qquad k\in\mathbb{Z},
$$
is given by the $3\times 3$ matrix~\footnote{Here and in what follows, we describe morphisms between finite biproducts in terms of matrices of morphisms such that the $j^{th}$-column of the matrix gives the morphisms from the $j^{th}$ factor of the domain to the various factors of the codomain.}
$$
d_k=\left(\begin{array}{ccc} 0&0&id_{B_k} \\ 0&0&0 \\ 0&0&0\end{array}\right). 
$$
Let us consider arbitrary morphisms in $\Aa$
$$
f_k:B_{k}\oplus H_{k}\oplus B_{k-1}\longrightarrow B_k\oplus H_k\oplus B_{k-1},\qquad k\in\mathbb{Z},
$$
described by matrices
$$
f_k=\left(\begin{array}{ccc} f(11)_k&f(12)_k&f(13)_k \\ f(21)_k&f(22)_k&f(23)_k \\ f(31)_k&f(32)_k&f(33)_k\end{array}\right).
$$

\begin{prop} \label{condicions_endomorfisme}
The sequence of morphisms $\{f_k\}_{k\in\mathbb{Z}}$ gives the components of an endomorphism $f$ of $A_{\bullet}$ if and only if
\begin{itemize}
\item $f(21)_k$, $f(31)_k$ and $f(32)_k$ are zero, and
\item $f(33)_k=f(11)_{k-1}$
\end{itemize}
for all $k\in\mathbb{Z}$. Moreover, the map
\begin{equation} \label{iso1}
End_{\Cc h(\Aa)}(A_{\bullet})\longrightarrow\prod_{k\in\mathbb{Z}}End_{\Aa}(B_k)\times End_{\Aa}(H_k)\times Hom_{\Aa}(H_k,B_k)\times Hom_{\Aa}(B_{k-1},H_k)\times Hom_{\Aa}(B_{k-1},B_k)
\end{equation}
given by $f\mapsto(f(11)_k,f(22)_k,f(12)_k,f(23)_k,f(13)_k)_{k\in\mathbb{Z}}$ is an isomorphism of $\kb$-modules.
\end{prop}
\begin{proof}
The endomorphism condition is $f_{k-1}\circ d_{k-1}=d_{k-1}\circ f_{k}$ for all $k\in\mathbb{Z}$. By taking the corresponding matrix products this gives
$$
\left(\begin{array}{ccc} 0&0&f(11)_{k-1} \\ 0&0&f(21)_{k-1} \\ 0&0&f(31)_{k-1}\end{array}\right)=\left(\begin{array}{ccc} f(31)_k&f(32)_k&f(33)_k \\ 0&0&0 \\ 0&0&0 \end{array}\right)\qquad\qquad \forall k\in\mathbb{Z},
$$
from which the first statement readily follows. Last assertion follows from the fact that there are no constraints on the remaining entries in $f_k$, and the fact that the sum and product by scalars between endomorphisms of $A_{\bullet}$ correspond to these same operations between the entries of the respective matrices. 
\end{proof}
From now on, we shall use the notation
\begin{equation} \label{matriu_f_k}
f_k=\left(\begin{array}{ccc} \phi_k&a_k&c_k \\ 0&\psi_k&b_k \\ 0&0&\phi_{k-1}\end{array}\right)
\end{equation}
for the matrices giving the components of an arbitrary endomorphism $f$ of $A_{\bullet}$. In particular, $\phi_k$ and $\psi_k$ are arbitrary endomorphisms of $B_k$ and $H_k$, respectively, while
\begin{align*}
a_k&:H_k\To B_k \\ b_k&:B_{k-1}\To H_k \\ c_k&:B_{k-1}\To B_k 
\end{align*}
are arbitrary morphisms. The image of $f$ by the isomorphism (\ref{iso1}) will be denoted by $(\mbox{\boldmath$\phi$\unboldmath},\mbox{\boldmath$\psi$\unboldmath},\mathbf{a},\mathbf{b},\mathbf{c})$. We can identify $f$ with its image by this isomorphism, and this is often done in what follows. In this case, the homotopy class of $f$ is denoted by $[\mbox{\boldmath$\phi$\unboldmath},\mbox{\boldmath$\psi$\unboldmath},\mathbf{a},\mathbf{b},\mathbf{c}]$. In particular, we have
\begin{equation} \label{id}
id_{A_\bullet}=(\mathbf{id},\mathbf{id},\mathbf{0},\mathbf{0},\mathbf{0}).
\end{equation}
Notice that the codomain of (\ref{iso1}) has a priori no $\kb$-algebra structure, but it gets one from the domain. The reader may easily check that the induced product is given by
\begin{equation} \label{producte}
(\mbox{\boldmath$\phi$\unboldmath}',\mbox{\boldmath$\psi$\unboldmath}',\mathbf{a}',\mathbf{b}',\mathbf{c}')\cdot(\mbox{\boldmath$\phi$\unboldmath},\mbox{\boldmath$\psi$\unboldmath},\mathbf{a},\mathbf{b},\mathbf{c})=(\mbox{\boldmath$\phi$\unboldmath}'\circ\mbox{\boldmath$\phi$\unboldmath},\mbox{\boldmath$\psi$\unboldmath}'\circ\mbox{\boldmath$\psi$\unboldmath},\mbox{\boldmath$\phi$\unboldmath}'\circ\mathbf{a}+\mathbf{a}'\circ\mbox{\boldmath$\psi$\unboldmath},\mbox{\boldmath$\psi$\unboldmath}'\circ\mathbf{b}+\mathbf{b}'\circ\mbox{\boldmath$\phi$\unboldmath},\mbox{\boldmath$\phi$\unboldmath}'\circ\mathbf{c}+\mathbf{a}'\circ\mathbf{b}+\mathbf{c}'\circ\mbox{\boldmath$\phi$\unboldmath})
\end{equation}
where
\begin{align*}
(\mbox{\boldmath$\phi$\unboldmath}'\circ\mbox{\boldmath$\phi$\unboldmath},&\mbox{\boldmath$\psi$\unboldmath}'\circ\mbox{\boldmath$\psi$\unboldmath},\mbox{\boldmath$\phi$\unboldmath}'\circ\mathbf{a}+\mathbf{a}'\circ\mbox{\boldmath$\psi$\unboldmath},\mbox{\boldmath$\psi$\unboldmath}'\circ\mathbf{b}+\mathbf{b}'\circ\mbox{\boldmath$\phi$\unboldmath},\mbox{\boldmath$\phi$\unboldmath}'\circ\mathbf{c}+\mathbf{a}'\circ\mathbf{b}+\mathbf{c}'\circ\mbox{\boldmath$\phi$\unboldmath})_k \\ &=(\phi'_k\circ\phi_k,\psi'_k\circ\psi_k,\phi'_k\circ a_k+a'_k\circ\psi_k,\psi'_k\circ b_k+b'_k\circ\phi_{k-1},\phi'_k\circ c_k+a'_k\circ b_k+c'_k\circ\phi_{k-1}).
\end{align*}
Although we shall not need it, let us remark that this formula allows us to identify the group of automorphisms of $A_{\bullet}$. More precisely, we have the following.

\begin{prop}
An endomorphism $f=(\mbox{\boldmath$\phi$\unboldmath},\mbox{\boldmath$\psi$\unboldmath},\mathbf{a},\mathbf{b},\mathbf{c})$ of $A_\bullet$ is an automorphism if and only if $\mbox{\boldmath$\phi$\unboldmath}$ and $\mbox{\boldmath$\psi$\unboldmath}$ are automorphisms of the $\mathbb{Z}$-graded objects $B_\bullet$ and $H_\bullet$, respectively (i.e. $\phi_k\in Aut_{\Aa}(B_k)$ and $\psi_k\in Aut_{\Aa}(H_k)$ for all $k\in\mathbb{Z}$).
\end{prop}
\begin{proof}
It follows from (\ref{producte}) that $(\mbox{\boldmath$\phi$\unboldmath},\mbox{\boldmath$\psi$\unboldmath},\mathbf{a},\mathbf{b},\mathbf{c})$ is an automorphism of $A_{\bullet}$ if and only if there exists $(\mbox{\boldmath$\phi$\unboldmath}',\mbox{\boldmath$\psi$\unboldmath}',\mathbf{a}',\mathbf{b}',\mathbf{c}')$ such that
\begin{eqnarray*}
&\phi'_k\circ\phi_k=id_{B_k} \\ &\psi'_k\circ\psi_k=id_{H_k} \\ &\phi'_k\circ a_k+a'_k\circ\psi_k=0 \\ &\psi'_k\circ b_k+b'_k\circ\phi_{k-1}=0 \\ &\phi'_k\circ c_k+a'_k\circ b_k+c'_k\circ\phi_{k-1}=0
\end{eqnarray*}  
for all $k\in\mathbb{Z}$. But this holds if and only if $\phi_k\in Aut_{\Aa}(B_k)$ and $\psi_k\in Aut_{\Aa}(H_k)$ for all $k\in\mathbb{Z}$. Indeed, the conditions are clearly necessary for the first two conditions to be satisfied, and they are also sufficient because the remaining three conditions automatically hold by taking
\begin{align*}
a'_k&=-\phi^{-1}_k\circ a_k\circ\psi^{-1}_k \\ b'_k&=-\psi^{-1}_k\circ b_k\circ\phi^{-1}_{k-1}
\end{align*}
and
$$
 c'_k=-\phi^{-1}_k\circ c_k\circ\phi^{-1}_{k-1}-a'_k\circ b_k\circ\phi^{-1}_{k-1}=-\phi^{-1}_k\circ c_k\circ\phi^{-1}_{k-1}+\phi^{-1}_k\circ a_k\circ\psi^{-1}_k\circ b_k\circ\phi^{-1}_{k-1}.
$$
\end{proof}
Therefore the underlying set of the group $\mathrm{Aut}_{\Cc h(\Aa)}(A_{\bullet})$ is
$$
\mathrm{Aut}_{\Cc h(\Aa)}(A_{\bullet})\cong\prod_{k\in\mathbb{Z}}Aut_{\Aa}(B_k)\times Aut_\Aa(H_k)\times Hom_\Aa(H_k,B_k)\times Hom_\Aa(B_{k-1},H_k)\times Hom_\Aa(B_{k-1},B_k).
$$
Observe that the group structure on this set is given by (\ref{producte}) and consequently, by the group structures of $Aut_\Aa(B_k)$, $Aut_\Aa(H_k)$ and $Aut_\Aa(B_{k-1})$ together with the canonical bimodule structures on $Hom_\Aa(H_k,B_k)$, $Hom_\Aa(B_{k-1},H_k)$ and $Hom_\Aa(B_{k-1},B_k)$.

Let us now determine when two endomorphisms of $A_{\bullet}$ are homotopic. Let $f=(\mbox{\boldmath$\phi$\unboldmath},\mbox{\boldmath$\psi$\unboldmath},\mathbf{a},\mathbf{b},\mathbf{c})$ and $f'=(\mbox{\boldmath$\phi'$\unboldmath},\mbox{\boldmath$\psi'$\unboldmath},\mathbf{a'},\mathbf{b'},\mathbf{c'})$. Then $f$ and $f'$ are homotopic if there exist morphisms
$$
h_k:A_k\longrightarrow A_{k+1},\qquad k\in\mathbb{Z},
$$
such that
$$
f'_k=f_k+d_k\circ h_k+h_{k-1}\circ d_{k-1},\qquad k\in\mathbb{Z}.
$$
In terms of matrices, this amounts to the existence of matrices of morphisms in $\Aa$
$$
h_k=\left(\begin{array}{ccc} h(11)_k&h(12)_k&h(13)_k \\ h(21)_k&h(22)_k&h(23)_k \\ h(31)_k&h(32)_k&h(33)_k\end{array}\right),\qquad k\in\mathbb{Z} 
$$
such that
$$
\left(\begin{array}{ccc} \phi'_k&a'_k&c'_k \\ 0&\psi'_k&b'_k \\ 0&0&\phi'_{k-1}\end{array}\right)=\left(\begin{array}{ccc} \phi_k&a_k&c_k \\ 0&\psi_k&b_k \\ 0&0&\phi_{k-1}\end{array}\right)+\left(\begin{array}{ccc} h(31)_k&h(32)_k&h(33)_k+h(11)_{k-1} \\ 0&0&h(21)_{k-1} \\ 0&0&h(31)_{k-1}\end{array}\right)
$$
and hence, such that
\begin{align*}
\phi'_k&=\phi_k+h(31)_k \\ a'_k&=a_k+h(32)_k \\ c'_k&=c_k+h(33)_k+h(11)_{k-1} \\ \psi'_k&=\psi_k \\ b'_k&=b_k+h(21)_{k-1}
\end{align*}
for all $k\in\mathbb{Z}$. Therefore we have the following:
\begin{prop} \label{relacio_homotopia}
Let $f=(\mbox{\boldmath$\phi$\unboldmath},\mbox{\boldmath$\psi$\unboldmath},\mathbf{a},\mathbf{b},\mathbf{c})$ and $f'=(\mbox{\boldmath$\phi'$\unboldmath},\mbox{\boldmath$\psi'$\unboldmath},\mathbf{a'},\mathbf{b'},\mathbf{c'})$ be arbitrary endomorphisms of $A_{\bullet}$. Then
$$
f\simeq f'\ \Longleftrightarrow\ \mbox{\boldmath$\psi$\unboldmath}=\mbox{\boldmath$\psi'$\unboldmath}.
$$
In this case, a homotopy is given by any collection of matrices 
$$
h_k=\left(\begin{array}{ccc} \alpha_k&\gamma_k&\delta_k \\ b'_{k+1}-b_{k+1}&\beta_k&\varepsilon_k \\ \phi'_k-\phi_k&a'_k-a_k&c'_k-c_k-\alpha_{k-1}\end{array}\right),\qquad k\in\mathbb{Z}, 
$$
with
\begin{align*}
\alpha_k&:B_k\longrightarrow B_{k+1} \\ \beta_k&:H_k\longrightarrow H_{k+1} \\ \gamma_k&:H_k\longrightarrow B_{k+1} \\ \delta_k&:B_{k-1}\longrightarrow B_{k+1} \\ \varepsilon_k&:B_{k-1}\longrightarrow H_{k+1}
\end{align*}
arbitrary morphisms in $\Aa$. 
\end{prop}
The above homotopy will be denoted by $h=(f,(\mbox{\boldmath$\alpha$\unboldmath},\mbox{\boldmath$\beta$\unboldmath},\mbox{\boldmath$\gamma$\unboldmath},\mbox{\boldmath$\delta$\unboldmath},\mbox{\boldmath$\varepsilon$\unboldmath}),f')$. Observe that the notation would be ambiguous without making explicit the domain and codomain of $h$.

\begin{cor}
An endomorphism $(\mbox{\boldmath$\phi$\unboldmath},\mbox{\boldmath$\psi$\unboldmath},\mathbf{a},\mathbf{b},\mathbf{c})$ of $A_\bullet$ is an equivalence  if and only if $\mbox{\boldmath$\psi$\unboldmath}$ is invertible, and in this case a pseudoinverse is given by the automorphism
\begin{equation} \label{pseudoinvers}
(\mbox{\boldmath$\phi$\unboldmath},\mbox{\boldmath$\psi$\unboldmath},\mathbf{a},\mathbf{b},\mathbf{c})^*=(\mathbf{id},\mbox{\boldmath$\psi$\unboldmath}^{-1},\mathbf{0},\mathbf{0},\mathbf{0})
\end{equation}
and also by the non strictly invertible morphism
\begin{equation} \label{pseudoinvers_bis}
(\mbox{\boldmath$\phi$\unboldmath},\mbox{\boldmath$\psi$\unboldmath},\mathbf{a},\mathbf{b},\mathbf{c})^*=(\mathbf{0},\mbox{\boldmath$\psi$\unboldmath}^{-1},\mathbf{0},\mathbf{0},\mathbf{0}).
\end{equation}
In particular, any self-equivalence of $A_{\bullet}$ is homotopic to an automorphism.
\end{cor}
\begin{proof}
By definition, $(\mbox{\boldmath$\phi$\unboldmath},\mbox{\boldmath$\psi$\unboldmath},\mathbf{a},\mathbf{b},\mathbf{c})$ is an equivalence if and only if there exists $(\mbox{\boldmath$\phi$\unboldmath}',\mbox{\boldmath$\psi$\unboldmath}',\mathbf{a}',\mathbf{b}',\mathbf{c}')$ such that
$$
(\mbox{\boldmath$\phi$\unboldmath}'\circ\mbox{\boldmath$\phi$\unboldmath},\mbox{\boldmath$\psi$\unboldmath}'\circ\mbox{\boldmath$\psi$\unboldmath},\mbox{\boldmath$\phi$\unboldmath}'\circ\mathbf{a}+\mathbf{a}'\circ\mbox{\boldmath$\psi$\unboldmath},\mbox{\boldmath$\psi$\unboldmath}'\circ\mathbf{b}+\mathbf{b}'\circ\mbox{\boldmath$\phi$\unboldmath},\mbox{\boldmath$\phi$\unboldmath}'\circ\mathbf{c}+\mathbf{a}'\circ\mathbf{b}+\mathbf{c}'\circ\mbox{\boldmath$\phi$\unboldmath})\simeq(\mathbf{id},\mathbf{id},\mathbf{0},\mathbf{0},\mathbf{0})
$$
$$
(\mbox{\boldmath$\phi$\unboldmath}\circ\mbox{\boldmath$\phi$\unboldmath}',\mbox{\boldmath$\psi$\unboldmath}\circ\mbox{\boldmath$\psi$\unboldmath}',\mbox{\boldmath$\phi$\unboldmath}\circ\mathbf{a}'+\mathbf{a}\circ\mbox{\boldmath$\psi$\unboldmath}',\mbox{\boldmath$\psi$\unboldmath}\circ\mathbf{b}'+\mathbf{b}\circ\mbox{\boldmath$\phi$\unboldmath}',\mbox{\boldmath$\phi$\unboldmath}\circ\mathbf{c}'+\mathbf{a}\circ\mathbf{b}'+\mathbf{c}\circ\mbox{\boldmath$\phi$\unboldmath}')\simeq(\mathbf{id},\mathbf{id},\mathbf{0},\mathbf{0},\mathbf{0}).
$$
The first statement follows now from the previous Proposition. As for the second statement, it follows from the previous Proposition and the above characterization of the automorphisms of $A_\bullet$.
\end{proof}
Later on we shall also need the following.

\begin{prop} \label{homotopies_homotopes_bis}
Let be given arbitrary endomorphisms $f=(\mbox{\boldmath$\phi$\unboldmath},\mbox{\boldmath$\psi$\unboldmath},\mathbf{a},\mathbf{b},\mathbf{c})$ and $f'=(\mbox{\boldmath$\phi'$\unboldmath},\mbox{\boldmath$\psi'$\unboldmath},\mathbf{a'},\mathbf{b'},\mathbf{c'})$, and let $h=(f,(\mbox{\boldmath$\alpha$\unboldmath},\mbox{\boldmath$\beta$\unboldmath},\mbox{\boldmath$\gamma$\unboldmath},\mbox{\boldmath$\delta$\unboldmath},\mbox{\boldmath$\varepsilon$\unboldmath}),f')$ and $h'=(f,(\mbox{\boldmath$\alpha'$\unboldmath},\mbox{\boldmath$\beta'$\unboldmath},\mbox{\boldmath$\gamma'$\unboldmath},\mbox{\boldmath$\delta'$\unboldmath},\mbox{\boldmath$\varepsilon'$\unboldmath}),f')$ be two homotopies between them. Then
$$
h\simeq h'\ \Leftrightarrow\ \mbox{\boldmath$\beta$\unboldmath}=\mbox{\boldmath$\beta'$\unboldmath}.
$$
In this case, a degree 2 homotopy between them is given by any collection of matrices
$$
h^{(2)}_k=\left(\begin{array}{ccc} \lambda_k&\nu_k&\theta_k \\ \varepsilon_{k+1}-\varepsilon'_{k+1}&\mu_k&\zeta_k \\ \alpha'_k-\alpha_k&\gamma'_k-\gamma_k&\delta'_k-\delta_k+\lambda_{k-1}\end{array}\right),\qquad k\in\mathbb{Z}, 
$$
with
\begin{align*}
\lambda_k&:B_k\longrightarrow B_{k+2} \\ \mu_k&:H_k\longrightarrow H_{k+2} \\ \nu_k&:H_k\longrightarrow B_{k+2} \\ \theta_k&:B_{k-1}\longrightarrow B_{k+2} \\ \zeta_k&:B_{k-1}\longrightarrow H_{k+2}
\end{align*}
arbitrary morphisms of $\Aa$. 
\end{prop}
\begin{proof}
By definition, $h$ and $h'$ are homotopic if there exists morphisms
$$
h^{(2)}_k:A_k\longrightarrow A_{k+2},\qquad k\in\mathbb{Z}
$$
such that
\begin{equation} \label{homotopies_homotopes}
h'_k=h_k+d_{k+1}\circ h^{(2)}_k-h^{(2)}_{k-1}\circ d_{k-1},\qquad k\in\mathbb{Z}.
\end{equation}
In terms of matrices, this amounts to the existence of matrices of morphisms in $\Aa$
$$
h^{(2)}_k=\left(\begin{array}{ccc} h^{(2)}(11)_k&h^{(2)}(12)_k&h^{(2)}(13)_k \\ h^{(2)}(21)_k&h^{(2)}(22)_k&h^{(2)}(23)_k \\ h^{(2)}(31)_k&h^{(2)}(32)_k&h^{(2)}(33)_k\end{array}\right),\qquad k\in\mathbb{Z} 
$$
such that
\begin{align*}
\alpha'_k&=\alpha_k+h^{(2)}(31)_k \\ \gamma'_k&=\gamma_k+h^{(2)}(32)_k \\ \delta'_k&=\delta_k+h^{(2)}(33)_k-h^{(2)}(11)_{k-1} \\ \beta'_k&=\beta_k \\ \varepsilon'_k&=\varepsilon_k-h^{(2)}(21)_{k-1} \\ \alpha'_{k-1}&=\alpha_{k-1}+h^{(2)}(31)_{k-1}
\end{align*}
for all $k\in\mathbb{Z}$. Notice that the last condition for all $k$ is equivalent to the first one and hence, redundant. Furthermore, the entries (2,1), (3,1), (3,2) and (3,3) of $h^{(2)}$ are uniquely determined by $h$ and $h'$. The remaining ones, however, can be chosen arbitrarily, and the above statement follows.
\end{proof}
The homotopy class of the homotopy $(f,(\mbox{\boldmath$\alpha$\unboldmath},\mbox{\boldmath$\beta$\unboldmath},\mbox{\boldmath$\gamma$\unboldmath},\mbox{\boldmath$\delta$\unboldmath},\mbox{\boldmath$\varepsilon$\unboldmath}),f')$ will be denoted by $[f,(\mbox{\boldmath$\alpha$\unboldmath},\mbox{\boldmath$\beta$\unboldmath},\mbox{\boldmath$\gamma$\unboldmath},\mbox{\boldmath$\delta$\unboldmath},\mbox{\boldmath$\varepsilon$\unboldmath}),f']$. In particular, the identity morphism of $f$ is 
$$
id_f=[f,(\mathbf{0},\mathbf{0},\mathbf{0},\mathbf{0},\mathbf{0}),f].
$$ 
We can now easily identify the $\kb$-algebra $H_0(End(V_{\bullet}))$ and its group of units.
\begin{prop} \label{pi_0}
Let $A_{\bullet}$ be the split chain complex defined by the objects $\{B_k,k\in\mathbb{Z}\}$ and $\{H_k,k\in\mathbb{Z}\}$ of $\Aa$. Then the map
\begin{equation} \label{iso_Psi}
\Psi:H_0(End(A_{\bullet}))\longrightarrow\prod_{k\in\mathbb{Z}}End_\Aa(H_k)
\end{equation}
given by $[(\mbox{\boldmath$\phi$\unboldmath},\mbox{\boldmath$\psi$\unboldmath},\mbox{\boldmath$a$\unboldmath},\mbox{\boldmath$b$\unboldmath},\mbox{\boldmath$c$\unboldmath})]\mapsto\mbox{\boldmath$\psi$\unboldmath}$ is an isomorphism of $\kb$-algebras. In particular, we have
\begin{equation} \label{iso_pi_0_G}
\pi_0(A_{\bullet})\cong\prod_{k\in\mathbb{Z}}Aut_\Aa(H_k).
\end{equation}
\end{prop}
\begin{proof}
It follows from Proposition~\ref{relacio_homotopia} that $\Psi$ is well defined injective map. Surjectivity follows from Proposition~\ref{condicions_endomorfisme}, which ensures that $\mbox{\boldmath$\psi$\unboldmath}$ can be chosen arbitrarily. Moreover, $\Psi$ is clearly linear and preserves the products as a consequence of (\ref{producte}). Last assertion follows because $\pi_0(A_{\bullet})$ is the group of units of $H_0(End(A_{\bullet}))$.
\end{proof}

The map (\ref{iso_Psi}) is nothing but the morphism $\Phi$ given by (\ref{morfisme_canonic}). Thus if $f=(\mbox{\boldmath$\phi$\unboldmath},\mbox{\boldmath$\psi$\unboldmath},\mathbf{a},\mathbf{b},\mathbf{c})$, it is easy to check that $H_k(f)=\psi_k$ for all $k\in\mathbb{Z}$. Hence, as claimed before, $\Phi$ is indeed an isomorphism of $\kb$-algebras when $A_{\bullet}=B_{\bullet}$ is a split chain complex.

\paragraph{\em The group $\pi_1(A_{\bullet})$.}

According to (\ref{pi_1_morfismes_A_A1}), $\pi_1(A_{\bullet})$ is the 0-homology of the complex $Hom(A_{\bullet},A[1]_{\bullet})$. To compute it, we proceed as before.

By definition, the boundary operator of $A[1]_{\bullet}$
$$
d[1]_{k-1}:B_{k+1}\oplus H_{k+1}\oplus B_k\longrightarrow B_k\oplus H_k\oplus B_{k-1}
$$
is given by
$$
d[1]_{k-1}=-d_k=-\left(\begin{array}{ccc} 0&0&id_{B_k} \\ 0&0&0 \\ 0&0&0\end{array}\right). 
$$
Let $g_k:A_k\To A[1]_k$ be arbitrary morphisms
$$
g_{k}:B_{k}\oplus H_{k}\oplus B_{k-1}\longrightarrow B_{k+1}\oplus H_{k+1}\oplus B_{k},\qquad k\in\mathbb{Z},
$$
described by the matrices
$$
g_k=\left(\begin{array}{ccc} g(11)_k&g(12)_k&g(13)_k \\ g(21)_k&g(22)_k&g(23)_k \\ g(31)_k&g(32)_k&g(33)_k\end{array}\right),\qquad k\in\mathbb{Z}.
$$
The analog of Lema~\ref{condicions_endomorfisme} reads now as follows:

\begin{prop} \label{condicions_morfisme_V_V1}
The sequence of morphisms $\{g_k\}_{k\in\mathbb{Z}}$ gives the components of a morphism $g:A_{\bullet}\To A[1]_{\bullet}$ if and only if
\begin{itemize}
\item $g(21)_k$, $g(31)_k$ and $g(32)_k$ are zero, and
\item $g(33)_k=-g(11)_{k-1}$
\end{itemize}
for all $k\in\mathbb{Z}$. Moreover, the map given by $g\mapsto\{(g(11)_k,g(22)_k,g(12)_k,g(23)_k,g(13)_k)\}_{k\in\mathbb{Z}}$ defines an isomorphism of $\kb$-modules
\begin{align*}
H&om_{\Cc h(\Aa)}(A_{\bullet},A[1]_{\bullet}) \\ &\cong\prod_{k\in\mathbb{Z}}Hom_\Aa(B_k,B_{k+1})\times Hom_\Aa(H_k,H_{k+1})\times Hom_\Aa(H_k,B_{k+1})\times Hom_\Aa(B_{k-1},H_{k+1})\times Hom_\Aa(B_{k-1},B_{k+1}).
\end{align*}
\end{prop}
We shall write (notice the minus sign in the (3,3)-component)
\begin{equation} \label{matriu_f'_k}
g_k=\left(\begin{array}{ccc} \rho_k&u_k&w_k \\ 0&\xi_k&v_k \\ 0&0&-\rho_{k-1}\end{array}\right)
\end{equation}
for any morphism $g:A_{\bullet}\To A[1]_{\bullet}$, with
\begin{align*}
\rho_k&:B_k\To B_{k+1} \\ \xi_k&:H_k\To H_{k+1} \\ u_k&:H_k\To B_{k+1} \\ v_k&:B_{k-1}\To H_{k+1} \\ w_k&:B_{k-1}\To B_{k+1} 
\end{align*}
arbitrary morphisms in $\Aa$. We shall denote by $(\mbox{\boldmath$\rho$\unboldmath},\mbox{\boldmath$\xi$\unboldmath},\mathbf{u},\mathbf{v},\mathbf{w})$ its image by the previous isomorphism, and we shall often identify $g$ with this image. Recall that we can also think of $g$ as a hotomopy of $id_{A_\bullet}$ to itself, in which case we should write $g=(id_{A_\bullet},(\mbox{\boldmath$\rho$\unboldmath},\mbox{\boldmath$\xi$\unboldmath},\mathbf{u},\mathbf{v},\mathbf{w}),id_{A_\bullet})$.

\begin{prop} \label{relacio_homotopia_bis}
Let $g=(\mbox{\boldmath$\rho$\unboldmath},\mbox{\boldmath$\xi$\unboldmath},\mathbf{u},\mathbf{v},\mathbf{w})$ and $g'=(\mbox{\boldmath$\rho'$\unboldmath},\mbox{\boldmath$\xi'$\unboldmath},\mathbf{u'},\mathbf{v'},\mathbf{w'})$ be arbitrary morphisms between $A_{\bullet}$ and $A[1]_{\bullet}$. Then 
$$
g\simeq g'\ \Longleftrightarrow\ \mbox{\boldmath$\xi$\unboldmath}=\mbox{\boldmath$\xi'$\unboldmath}.
$$
In this case, a homotopy is given by any collection of matrices
$$
h_k=\left(\begin{array}{ccc} \lambda_k&\nu_k&\theta_k \\ v'_{k+1}-v_{k+1}&\mu_k&\zeta_k \\ \rho_k-\rho'_k&u_k-u'_k&w_k-w'_k+\lambda_{k-1}\end{array}\right),\qquad k\in\mathbb{Z}, 
$$
with
\begin{align*}
\lambda_k&:B_k\longrightarrow B_{k+2} \\ \mu_k&:H_k\longrightarrow H_{k+2} \\ \nu_k&:H_k\longrightarrow B_{k+2} \\ \theta_k&:B_{k-1}\longrightarrow B_{k+2} \\ \zeta_k&:B_{k-1}\longrightarrow H_{k+2}
\end{align*}
arbitrary morphisms of $\Aa$. 
\end{prop}
\begin{proof}
The proof is the same as that of Proposition~\ref{homotopies_homotopes_bis}, except that instead of (\ref{homotopies_homotopes}) we now have the condition
$$
g'_k=g_k+d[1]_k\circ h_k+h_{k-1}\circ d_{k-1},\qquad k\in\mathbb{Z}.
$$
The details are left to the reader.
\end{proof}
The homotopy class of $(\mbox{\boldmath$\rho$\unboldmath},\mbox{\boldmath$\xi$\unboldmath},\mathbf{u},\mathbf{v},\mathbf{w})$ will be denoted by $[\mbox{\boldmath$\rho$\unboldmath},\mbox{\boldmath$\xi$\unboldmath},\mathbf{u},\mathbf{v},\mathbf{w}]$.
 
We can now identify the abelian group $\pi_1(A_{\bullet})$.
\begin{prop}
Let $A_{\bullet}$ be the split chain complex defined by the objects $\{B_k,k\in\mathbb{Z}\}$ and $\{H_k,k\in\mathbb{Z}\}$ of $\Aa$. Then the map
$$
\Theta:H_0(Hom(A_{\bullet},A[1]_{\bullet}))\longrightarrow\prod_{k\in\mathbb{Z}}Hom_\Aa(H_k,H_{k+1})
$$
given by $[\mbox{\boldmath$\rho$\unboldmath},\mbox{\boldmath$\xi$\unboldmath},\mbox{\boldmath$u$\unboldmath},\mbox{\boldmath$v$\unboldmath},\mbox{\boldmath$w$\unboldmath}]\mapsto\mbox{\boldmath$\xi$\unboldmath}$ is an isomorphism de $\kb$-modules. In particular we have
\begin{equation} \label{iso_pi_1_G}
\pi_1(A_{\bullet})\cong\prod_{k\in\mathbb{Z}}Hom_\Aa(H_k,H_{k+1}).
\end{equation}
\end{prop}

\begin{rem}{\rm
More generally, for any $n\geq 1$ it can be shown that $H_0(Hom(A_\bullet,A[n]_\bullet))$ is isomorphic to $\prod_{k\in\mathbb{Z}}Hom_\Aa(H_k,H_{k+n})$. }
\end{rem}
Combined with Proposition~\ref{pi_0}, this result already proves that, up to equivalence, the 2-group of symmetries of any split exact sequence is trivial. Hence the non-triviality of the 2-group of symmetries of a split chain complex can be thought of as a measure of its non-exactness.

\paragraph{\em Structure of $\pi_0(A_{\bullet})$-module on $\pi_1(A_{\bullet})$.}
Let us identify $\pi_0(A_{\bullet})$ and $\pi_1(A_{\bullet})$ with the above product groups (\ref{iso_pi_0_G}) and (\ref{iso_pi_1_G}). Under these identifications, (\ref{accio_pi0_pi1}) translates into the action
$$
\mbox{\boldmath$\psi$\unboldmath}\lhd\mbox{\boldmath$\xi$\unboldmath}=\mbox{\boldmath$\psi$\unboldmath}^{-1}[1]\circ\mbox{\boldmath$\xi$\unboldmath}\circ\mbox{\boldmath$\psi$\unboldmath}
$$
where
$$
(\mbox{\boldmath$\psi$\unboldmath}^{-1}[1]\circ\mbox{\boldmath$\xi$\unboldmath}\circ\mbox{\boldmath$\psi$\unboldmath})_k=\psi^{-1}_{k+1}\circ\xi_k\circ\psi_k
$$
for all $k\in\mathbb{Z}$.

\paragraph{\em Postnikov invariant.}
It follows from (\ref{producte}) that the maps $s_0,s_1:\pi_0(A_\bullet)\To Equiv(A_\bullet)$ defined by
$$
s_1([\mbox{\boldmath$\phi$\unboldmath},\mbox{\boldmath$\psi$\unboldmath},\mathbf{a},\mathbf{b},\mathbf{c}])=(\mathbf{id},\mbox{\boldmath$\psi$\unboldmath},\mathbf{0},\mathbf{0},\mathbf{0})
$$
$$
s_0([\mbox{\boldmath$\phi$\unboldmath},\mbox{\boldmath$\psi$\unboldmath},\mathbf{a},\mathbf{b},\mathbf{c}])=(\mathbf{0},\mbox{\boldmath$\psi$\unboldmath},\mathbf{0},\mathbf{0},\mathbf{0}).
$$
are both sections of the canonical projection $\pi:Equiv(A_{\bullet})\To\pi_0(A_{\bullet})$ in the category of monoids. $s_1$ is a section by automorphisms of $A_{\bullet}$ while $s_0$ is a section by non-strictly invertible self-equivalences. It follows from the general discussion in \S~\ref{cas_general} that the Postnikov invariant of $\mathbb{E}quiv(A_{\bullet})$ is zero for any split chain complex $A_{\bullet}$.

\paragraph{\em An explicit equivalence.}
Let us keep identifying $\pi_0(A_{\bullet})$ and $\pi_1(A_{\bullet})$ with the above product groups (\ref{iso_pi_0_G}) and (\ref{iso_pi_1_G}). To obtain an explicit equivalence of 2-groups
\begin{equation} \label{equivalencia}
(F,\mu):Hom_{\Aa}(H_\bullet,H_{\bullet+1})[1]\rtimes Aut_{\Aa}(H_\bullet)[0]\stackrel{\simeq}{\longrightarrow}\mathbb{E}quiv(A_\bullet)
\end{equation}
we follow the discussion in \S~\ref{classificacio}. First of all we need to choose representative objects for the elements in $\pi_0(A_\bullet)$ and an isomorphism between any self-equivalence of $A_\bullet$ and the chosen representative in its homotopy class. We shall freely use the notations introduced previously.

For any $\mbox{\boldmath$\psi$\unboldmath}\in Aut_\Aa(H_\bullet)$ (a homotopy class of self-equivalences of $A_\bullet$) we choose as its representative the self-equivalence
\begin{equation} \label{f_psi}
f_{\mbox{\boldmath$\psi$\unboldmath}}=(\mathbf{id},\mbox{\boldmath$\psi$\unboldmath},\mathbf{0},\mathbf{0},\mathbf{0}).
\end{equation}
In particular, we have
$$
f_{\mathbf{id}}=(\mathbf{id},\mathbf{id},\mathbf{0},\mathbf{0},\mathbf{0})=id_{A_\bullet}
$$
as required. For any other self-equivalence $f\in\mbox{\boldmath$\psi$\unboldmath}$ let
$$
\iota_f:f\To f_{\mbox{\boldmath$\psi$\unboldmath}}
$$
be the isomorphism given by the homotopy class of homotopies
\begin{equation} \label{iota_f}
\iota_f=[f,(\mathbf{0},\mathbf{0},\mathbf{0},\mathbf{0},\mathbf{0}),f_{\mbox{\boldmath$\psi$\unboldmath}}].
\end{equation}
Let us emphasize that $(f,(\mathbf{0},\mathbf{0},\mathbf{0},\mathbf{0},\mathbf{0}),f_{\mbox{\boldmath$\psi$\unboldmath}})$ is not the zero homotopy in general. In fact, the degree 1 zero homotopy does not define a morphism between $f$ and $f_{\mbox{\boldmath$\psi$\unboldmath}}$ unless $f=f_{\mbox{\boldmath$\psi$\unboldmath}}$. For an arbitrary $f=(\mbox{\boldmath$\phi$\unboldmath},\mbox{\boldmath$\psi$\unboldmath},\mathbf{a},\mathbf{b},\mathbf{c})$ in the same homotopy class as $f_{\mbox{\boldmath$\psi$\unboldmath}}$, $(f,(\mathbf{0},\mathbf{0},\mathbf{0},\mathbf{0},\mathbf{0}),f_{\mbox{\boldmath$\psi$\unboldmath}})$ stands for the degree 1 homotopy given by the matrices
$$
h_k=\left(\begin{array}{ccc} 0&0&0 \\ -b_{k+1}&0&0 \\ id_{B_k}-\phi_k&-a_k&-c_k\end{array}\right),\qquad k\in\mathbb{Z}
$$
(see Proposition~\ref{relacio_homotopia}). We choose $\iota_f$ as in (\ref{iota_f}) because we then have
$$
\iota_{f_{\mbox{\boldmath$\psi$\unboldmath}}}=[f_{\mbox{\boldmath$\psi$\unboldmath}},(\mathbf{0},\mathbf{0},\mathbf{0},\mathbf{0},\mathbf{0}),f_{\mbox{\boldmath$\psi$\unboldmath}}]=(f_{\mbox{\boldmath$\psi$\unboldmath}},[0])=id_{f_{\mbox{\boldmath$\psi$\unboldmath}}}
$$
as required.

According to (\ref{functor_F}), an equivalence (\ref{equivalencia}) is then given by the functor $F$ acting as follows. It maps the object $\mbox{\boldmath$\psi$\unboldmath}\in Aut_{\Aa}(H_\bullet)$ to $F(\mbox{\boldmath$\psi$\unboldmath})=f_{\mbox{\boldmath$\psi$\unboldmath}}$, i.e. the self-equivalence of $A_\bullet$
$$
\xymatrix{
\cdots\ar[rr] && B_k\oplus H_k\oplus B_{k-1}\ar[rr]\ar[dd]^{\left(\begin{array}{ccc} id_{B_k}&0&0 \\ 0&\psi_k&0 \\ 0&0&id_{B_{k-1}}\end{array}\right)} &&  \cdots \\ &&&& \\ \cdots\ar[rr] && B_k\oplus H_k\oplus B_{k-1}\ar[rr] && \cdots\ , }
$$
and the morphism $(\mbox{\boldmath$\xi$\unboldmath},\mbox{\boldmath$\psi$\unboldmath})\in Hom_{\Aa}(H_\bullet,H_{\bullet+1})\times Aut_{\Aa}(H_\bullet)$ to
$$
F(\mbox{\boldmath$\xi$\unboldmath},\mbox{\boldmath$\psi$\unboldmath})=\gamma_{f_{\mbox{\boldmath$\psi$\unboldmath}}}([\mbox{\boldmath$\rho$\unboldmath},\mbox{\boldmath$\xi$\unboldmath},\mathbf{u},\mathbf{v},\mathbf{w}]):f_{\mbox{\boldmath$\psi$\unboldmath}}\To f_{\mbox{\boldmath$\psi$\unboldmath}}.
$$
Here $\mbox{\boldmath$\rho$\unboldmath}$, $\mathbf{u}$, $\mathbf{v}$ and $\mathbf{w}$ can be chosen arbitrarily (cf. Proposition~\ref{relacio_homotopia_bis}). For the sake of simplicity, we shall take them equal to zero, and we shall denote by $h_{\mbox{\boldmath$\xi$\unboldmath}}$ the corresponding homotopy from $id_{A_\bullet}$ to itself. It is given by the matrices
$$
h_{\mbox{\boldmath$\xi$\unboldmath},k}=\left(\begin{array}{ccc} 0&0&0 \\ 0&\xi_k&0 \\ 0&0&0\end{array}\right),\qquad k\in\mathbb{Z}.
$$
Now, since the underlying monoidal groupoid of $\mathbb{E}quiv(A_\bullet)$ is strict we have
\begin{align*}
\gamma_{f_{\mbox{\boldmath$\psi$\unboldmath}}}([\mathbf{0},\mbox{\boldmath$\xi$\unboldmath},\mathbf{0},\mathbf{0},\mathbf{0}])\ &\stackrel{(\ref{gamma_x})}{=}\ id_{f_{\mbox{\boldmath$\psi$\unboldmath}}}\otimes(id_{A_\bullet},[h_{\mbox{\boldmath$\xi$\unboldmath}}]) \\ &\hspace{0.1 truecm} =\ (f_{\mbox{\boldmath$\psi$\unboldmath}},[0])\otimes(id_{A_\bullet},[h_{\mbox{\boldmath$\xi$\unboldmath}}]) \\ &\stackrel{(\ref{hat_h})}{=}\ (f_{\mbox{\boldmath$\psi$\unboldmath}},[h_{\mbox{\boldmath$\psi$\unboldmath}\circ\mbox{\boldmath$\xi$\unboldmath}}]),
\end{align*}
where $h_{\mbox{\boldmath$\psi$\unboldmath}\circ\mbox{\boldmath$\xi$\unboldmath}}$ is the degree 1 homotopy of $A_{\bullet}$ into itself given by the matrices
$$
h_{\mbox{\boldmath$\psi$\unboldmath}\circ\mbox{\boldmath$\xi$\unboldmath},k}=\left(\begin{array}{ccc} 0&0&0 \\ 0&\psi_{k+1}\circ\xi_k&0 \\ 0&0&0\end{array}\right),\qquad k\in\mathbb{Z}.
$$
In summary, $F(\mbox{\boldmath$\xi$\unboldmath},\mbox{\boldmath$\psi$\unboldmath})$ is the homotopy class of the degree 1 homotopy of $A_\bullet$ into itself
$$
\xymatrix{
&  \cdots\ar[r] & B_k\oplus H_k\oplus B_{k-1}\ar[r]\ar[dddl]_{\left(\begin{array}{ccc} 0&0&0 \\ 0&\psi_{k+1}\circ\xi_k&0 \\ 0&0&0\end{array}\right)} &  \cdots \\   &&& \\ &&& \\ \cdots\ar[r] & B_{k+1}\oplus H_{k+1}\oplus B_{k}\ar[r] & \cdots & & }
$$
(it is immediate to check that this is indeed an automorphism of $f_{\mbox{\boldmath$\psi$\unboldmath}}$). As for the monoidal structure of this functor, for any $\mbox{\boldmath$\psi$\unboldmath},\mbox{\boldmath$\psi'$\unboldmath}\in Aut_\Aa(H_\bullet)$ we have
$$
f_{\mbox{\boldmath$\psi$\unboldmath}}\circ f_{\mbox{\boldmath$\psi'$\unboldmath}}=(\mathbf{id},\mbox{\boldmath$\psi$\unboldmath},\mathbf{0},\mathbf{0},\mathbf{0})\circ(\mathbf{id},\mbox{\boldmath$\psi'$\unboldmath},\mathbf{0},\mathbf{0},\mathbf{0})=(\mathbf{id},\mbox{\boldmath$\psi$\unboldmath}\circ\mbox{\boldmath$\psi'$\unboldmath},\mathbf{0},\mathbf{0},\mathbf{0})=f_{\mbox{\boldmath$\psi$\unboldmath}\circ\mbox{\boldmath$\psi'$\unboldmath}}.
$$
Hence the isomorphism $\iota_{f_{\mbox{\boldmath$\psi$\unboldmath}}\circ f_{\mbox{\boldmath$\psi'$\unboldmath}}}$ is an identity and (\ref{mu}) then implies that $F$ is a strict monoidal functor.

\section{Appendix}

\subsection{Definition of a (weak) 2-category}

For the reader's convenience, we recall here the definition of a weak 2-category, also called a {\it bicategory} and often just a 2-category (see \cite{jB67}, \cite{KS74} or \cite{sM98}). Roughly, it is a category whose hom-sets are categories, and whose composition maps are functors. The weak condition corresponds to the fact that the associativity of composition and the unit character of the identity morphisms are assumed to hold only up to isomorphism. This leads to the following definition.

A (small) weak 2-category $\mathbf{C}$ consists consists of the following set of data:
\begin{itemize}
\item
A set $|\mathbf{C}|$ of objects.
\item
For any ordered pair of objects $X,Y\in|\mathbf{C}|$, a small category
$\Cgg(X,Y)$ whose objects (called {\it 1-morphisms}) are denoted by
$f:X\longrightarrow Y$ and whose morphisms (called {\it 2-morphisms}) are denoted by
$\tau:f\Longrightarrow f'$. The identity 2-morphism of a 1-morphism $f$ is
denoted by $1_f$ and the (strictly associative) composition between 2-morphisms (called
{\it vertical composition}) by $\tau'\cdot\tau$.
\item
For any ordered triple of objects $X,Y,Z\in|\Cgg|$, a functor
$$
{\bf comp}_{X,Y,Z}:\Cgg(X,Y)\times \Cgg(Y,Z)\longrightarrow
\Cgg(X,Z)
$$
whose action on objects $(f,g)$ is denoted by $g\circ f$ and whose action on
morphisms $(\tau,\sigma)$ is denoted by $\sigma\circ\tau$ and called
{\it horizontal composition}.
\item
For any object $X\in|\Cgg|$, a distinguished 1-morphism
$id_X\in|\Cgg(X,X)|$.
\item
For any objects $X,Y,Z,T\in|\Cgg|$ and any composable 1-morphisms
$f:X\longrightarrow Y$, $g:Y\longrightarrow Z$, $h:Z\longrightarrow
T$, a 2-isomorphism $\alpha_{h,g,f}:h\circ (g\circ
f)\Longrightarrow (h\circ g)\circ f$ (called {\it associativity constraint}) natural
in $f,g,h$.
\item
For any 1-morphism $f:X\longrightarrow Y$, two 2-isomorphisms
$\lambda_f:id_Y\circ f\Longrightarrow f$ and $\rho_f:f\circ
id_X\Longrightarrow f$ (respectively called {\it left} and {\it right unit
constraints}) natural in $f$.
\end{itemize}
Moreover, these data must be such that the diagrams (of vertical compositions)
$$
\xymatrix{
k\circ(h\circ(g\circ f))\ar@{=>}[d]_{\alpha_{k,h,g\circ f}}
\ar@{=>}[rr]^{1_k\circ \alpha_{h,g,f}} & & k\circ((h\circ g)\circ
f)\ar@{=>}[d]^{\alpha_{k,h\circ g,f}}
\\ (k\circ h)\circ (g\circ f)\ar@{=>}[dr]_{\alpha_{k\circ h,g,f}} &
& (k\circ(h\circ g))\circ f \ar@{=>}[dl]^{\alpha_{k,h,g}\circ 1_f}
\\ & ((k\circ h)\circ g)\circ f }
$$
$$
\xymatrix{
(g\circ id_Y)\circ f\ar@{=>}[rr]^{\alpha_{g,id_Y,f}}
\ar@{=>}[dr]_{\rho_g\circ 1_f} & & g\circ(id_Y\circ f)
\ar@{=>}[dl]^{1_g\circ \lambda_f}
\\ & g\circ f } 
$$
commute for all 1-morphisms $f,g,h$. When all isomorphism constraints $\alpha_{h,g,f}$, $\lambda_f$ and $\rho_f$ are identities we speak of a {\it strict} 2-category.

In this work we are only concerned with strict 2-categories of more than one object and arbitrary one-object 2-categories. These are the same thing as monoidal categories, the tensor product being given by the (unique) composition functor and the unit object by the identity 1-morphism of the unique object.

The fact that in a 2-category there are morphisms between morphisms implies that two objects $X,Y$ can be equal ($X=Y$), isomorphic ($X\cong Y$) or just {\it equivalent} ($X\simeq Y$), i.e. such that there exists 1-morphisms $f:X\To Y$ and $f^*:Y\To X$ such $f\circ f^*$ and $f^*\circ f$ are 2-isomorphic to the corresponding identity 1-morphisms.

\subsection{More on 2-groups}

\paragraph{\em On the canonical isomorphisms $\delta_x,\gamma_x$.}
As pointed out before, for any 2-group $\mathbb{G}$ we have an action of $\pi_0(\mathbb{G})$ on $\pi_1(\mathbb{G})$ defined by (cf. (\ref{accio_pi0}))
$$
[x]\lhd u=\gamma_x^{-1}(\delta_x(u)).
$$
Here $\gamma_x,\delta_x:\pi_1(\mathbb{G})\To Aut_{\Gg}(x)$ denote the canonical isomorphisms of groups induced by the monoidal structure on $\Gg$, and respectively given by
\begin{align*}
\delta_x(u)&=r_x\circ(u\otimes id_x)\circ r_x^{-1}, \\ \gamma_x(u)&=l_x\circ(id_x\otimes u)\circ l_x^{-1}. 
\end{align*}
An explicit expression for the corresponding inverse morphisms when $\mathbb{G}$ is such that the underlying monoidal groupoid $\Gg$ is strict is the following.  

\begin{prop}
Let $\mathbb{G}$ be a non necessarily strict 2-group whose underlying monoidal groupoid $\Gg$ is strict, and let $x$ be any object of $\Gg$. Then for any automorphism $\varphi:x\To x$ we have
\begin{align}
\gamma_x^{-1}(\varphi)&=\epsilon\circ(id_{x^*}\otimes\varphi)\circ\epsilon^{-1}, \label{gamma-1} \\ \delta_x^{-1}(\varphi)&=\eta^{-1}\circ(\varphi\otimes id_{x^*})\circ\eta, \label{delta-1}
\end{align}
where $x^*$ is any pseudoinverse of $x$ and $\epsilon$, $\eta$ are any isomorphism $\epsilon:x^*\otimes x\stackrel{\cong}{\To} I$ and $\eta:I\To x\otimes x^*$.
\end{prop}
\begin{proof} 
Given any triple $(x,x^*,\epsilon)$ as in the statement, it can be completed (in a unique way) to an adjoint equivalence $(x,x^*,\eta,\epsilon)$, i.e. there exists a (unique) isomorphism $\eta:I\To x\otimes x^*$ such that the composite morphisms
$$
\xymatrix{
x\ar[r]^{\cong\ \ \ \ } & I\otimes x\ar[r]^{\eta\otimes id_x\ \ \ \ } & (x\otimes x^*)\otimes x\ar[r]^{\cong} & x\otimes(x^*\otimes x)\ar[r]^{\ \ \ \ id_x\otimes\epsilon} & x\otimes I\ar[r]^{\cong} & x} 
$$
$$
\xymatrix{ x^*\ar[r]^{\cong\ \ \ \ } & x^*\otimes I\ar[r]^{id_{x^*}\otimes\eta\ \ \ \ } & x^*\otimes (x\otimes x^*)\ar[r]^{\cong} & (x^*\otimes x)\otimes x^*\ar[r]^{\ \ \ \ \epsilon\otimes id_{x^*}} & I\otimes x^*\ar[r]^{\cong} & x^*
}
$$
are both identities. In case $\Gg$ is strict, this means that
\begin{equation} \label{eps}
id_x\otimes\epsilon=(\eta\otimes id_x)^{-1}=\eta^{-1}\otimes id_x  
\end{equation}
because $\eta$ is invertible. Hence
\begin{align*}
\gamma_x(\epsilon\circ(id_{x^*}\otimes\varphi)\circ\epsilon^{-1})&\ =\ (id_x\otimes\epsilon)\circ(id_x\otimes id_{x^*}\otimes\varphi)\circ(id_x\otimes\epsilon^{-1}) \\ &\stackrel{(\ref{eps})}{=}(\eta^{-1}\otimes id_x)\circ(id_{x\otimes x^*}\otimes\varphi)\circ(id_x\otimes\epsilon^{-1}) \\ &\ =\ (\eta^{-1}\otimes\varphi)\circ(id_x\otimes\epsilon^{-1}) \\ &\ =\ (id_I\otimes\varphi)\circ(\eta^{-1}\otimes id_x)\circ(id_x\otimes\epsilon^{-1}) \\ &\ =\ \varphi\circ[(id_x\otimes\epsilon)\circ(\eta\otimes id_x)]^{-1} \\ &\ =\ \varphi.
\end{align*}
This proves (\ref{gamma-1}). (\ref{delta-1}) is shown in a similar way.
\end{proof}

\begin{cor}
Let $\mathbb{G}$ be a (non necessarily strict) 2-group whose underlying monoidal groupoid is strict. Then the action (\ref{accio_pi0}) of $\pi_0(\mathbb{G})$ on $\pi_1(\mathbb{G})$ is given by
$$
[x]\lhd u=\epsilon\circ(id_{x^*}\otimes u\otimes id_x)\circ\epsilon^{-1}
$$
for any representative $x$ of $[x]$, any pseudoinverse $x^*$ of $x$ and any isomorphism $\epsilon:x^*\otimes x\stackrel{\cong}{\To} I$.
\end{cor}

\paragraph{\em Sinh's theorem.}

Let us now prove that the pair $(F,\mu)$ defined by (\ref{functor_F})-(\ref{mu}) indeed defines an equivalence of 2-groups. We freely use the notations introduced in \S~\ref{classificacio}. Furthermore, for any object $x$ of $\Gg$ we shall denote by $g_x$ the corresponding isomorphism class, i.e. $g_x=[x]\in\pi_0(\mathbb{G})$. Notice that
$$
g_{x_g}=g
$$
for any $g\in\pi_0(\mathbb{G})$, whereas for an arbitrary object $x$ of $\Gg$ we have in general
$$
x_{g_x}\neq x
$$
because the chosen representative of $g_x$ need not be the object $x$. Equality holds if and only if $x$ is one of the chosen representatives. Thus
$$
x_{g_{x_g}}=x_g
$$
for all $g\in\pi_0(\mathbb{G})$.

Let $F^*:\Gg\To\Gg_{\pi_1(\mathbb{G}),\pi_0(\mathbb{G})}$ be the functor defined on objects $x$ and morphisms $\varphi:x\To y$ by
\begin{equation} \label{definicio_E'}
F^*(x)=g_x,\qquad F^*(\varphi)=(\gamma^{-1}_{x_{g_x}}(\iota_y\circ\varphi\circ\iota^{-1}_x),g_x).
\end{equation}
Observe that $\iota_y\circ\varphi\circ\iota^{-1}_x$ is a morphism from $x_{g_x}$ to $x_{g_y}$ and hence, an automorphism of $x_{g_x}$ because $g_x=g_y$.

\begin{lem}
$F^*\circ F=id$.
\end{lem}
\begin{proof}
For any object $g\in\pi_0(\mathbb{G})$ we have
$$
F^*(F(g))=F^*(x_g)=g_{x_g}=g,
$$
and for any morphism $(u,g):g\To g$ we have
$$
F^*(F(u,g))=F^*(\gamma_{x_g}^{-1}(u))=(\gamma_{x_{g_{x_g}}}^{-1}(\iota_{x_g}\circ\gamma_{x_g}^{-1}(u)\circ\iota_{x_g}^{-1}),g)=(\gamma_{x_g}^{-1}(\gamma_{x_g}(u)),g)=(u,g)
$$
because $x_{g_{x_g}}=x_g$ and $\iota_{x_g}$  is an identity.
\end{proof}

\begin{lem}
$F\circ F^*\cong id$.
\end{lem}
\begin{proof}
For any object $x$ of $\Gg$ we have
$$
F(F^*(x))=F(g_x)=x_{g_x}.
$$
Now, as pointed out before, in general $x_{g_x}$ is only isomorphic to $x$ through the isomorphism $\iota_x:x\To x_{g_x}$. We need to see that these are the components of a natural transformation $\iota:id\Rightarrow\Ee\circ\Ee'$. But this is an immediate consequence of the definitions. Thus for any morphism $\varphi:x\To y$ we have
$$
F(F^*(\varphi))=F(\gamma^{-1}_{x_{g_x}}(\iota_y\circ\varphi\circ\iota^{-1}_x),g_x)=\gamma_{x_{g_x}}(\gamma^{-1}_{x_{g_x}}(\iota_y\circ\varphi\circ\iota^{-1}_x))=\iota_y\circ\varphi\circ\iota^{-1}_x,
$$
and this is precisely the naturality of $\iota_x$ in $x$.
\end{proof}
This proves that $F$ is an equivalence of categories with $F^*$ as pseudoinverse. It remains to be shown that it is an equivalence of {\it monoidal} categories, for which it is enough to see that the natural isomorphism $\mu$ defined (\ref{mu}) indeed satisfies axiom (\ref{axioma_coherencia}). But this readily follows from (\ref{alpha}), (\ref{associador_G}) and (\ref{functor_F}).

\bibliographystyle{plain}
\bibliography{symmetries_chain_complex}

\begin{thebibliography}{10}

\bibitem{BC03}
J.~Baez and A.~Crans.
\newblock Higher-dimensional algebra {VI}: Lie 2-algebras.
\newblock {\em Theory Appl. Categ.}, 12:492--538, 2004 (also available as
  arXiv: math.QA/0307263).

\bibitem{BL03}
J.~Baez and A.~Lauda.
\newblock Higher-dimensional algebra {V}: 2-groups.
\newblock {\em Theory Appl. Categ.}, 12:423--491, 2004 (also available as
  arXiv: math.QA/0307200).

\bibitem{jB67}
J.~Benabou.
\newblock Introduction to bicategories.
\newblock In {\em Reports of the Midwest Category Seminar (LNM, volume 47)},
  pages 1--77. Springer, 1967.

\bibitem{jE4}
J.~Elgueta.
\newblock Representation theory of 2-groups on {K}apranov and {V}oevodsky
  2-vector spaces.
\newblock {\em Adv. Math.}, 213:53--92, 2007 (previous version available as
  arXiv.org: math.CT/0408120).

\bibitem{jE5}
J.~Elgueta.
\newblock Generalized 2-vector spaces and general linear 2-groups.
\newblock {\em J. Pure Appl. Alg.}, 212:2067--2091, 2008.

\bibitem{jE6}
J.~Elgueta.
\newblock On the regular representation of an (essentially) finite 2-group.
\newblock {\em Preprint {\rm arXiv:0907.0978}}, 2009.

\bibitem{GZ67}
P.~Gabriel and M.~Zisman.
\newblock {\em Calculus of fractions and homotopy theory}.
\newblock Springer Verlag, 1967.

\bibitem{GI01}
A.~Garzon and H.~Inassaridze.
\newblock Semidirect product of categorical groups, obstruction theory.
\newblock {\em Hom., Hom. and Applications}, 3:111--138, 2001.

\bibitem{KV94}
M.~Kapranov and V.~Voevodsky.
\newblock 2-categories and {Z}amolodchikov tetrahedra equations.
\newblock In {\em Proc. Sympos. Pure Math.}, volume 56(2), pages 177--260.
  American Mathematical Society, 1994.

\bibitem{KS74}
G.M. Kelly and R.~Street.
\newblock Review of the elements of 2-categories.
\newblock In {\em Proceedings of the Category Seminar, Sydney (LNM, volume
  420)}, pages 75--103. Springer, 1974.

\bibitem{LP08}
S.~Lack and S.~Paoli.
\newblock 2-nerves for bicategories.
\newblock {\em K-Theory}, 38:153--175, 2008.

\bibitem{sM98}
S.~MacLane.
\newblock {\em Categories for the Working Mathematician}, volume~5 of {\em
  Graduate Texts in Mathematics}.
\newblock Springer, 1998.

\bibitem{SR72}
N.~Saavedra Rivano.
\newblock {\em Cat\'{e}gories tannakiennes}, volume 265 of {\em Lecture Notes
  in Mathematics}.
\newblock Springer-Verlag, 1972.

\bibitem{hxSi75}
Hoang~Xuan Sinh.
\newblock {\em Gr-cat\'egories}.
\newblock Th\`ese de doctorat. Universit\'e Paris-VII, 1975.

\bibitem{cW94}
C.~Weibel.
\newblock {\em An introduction to homological algebra}, volume~28 of {\em
  Cambridge studies in advanced mathematics}.
\newblock Cambridge University Press, 1994.

\end{thebibliography}

\end{document}